\documentclass{article}
\usepackage{amssymb,amsbsy,amsmath}
\usepackage{natbib}
\usepackage{verbatim}
\usepackage{multirow}
\usepackage{graphicx}
\usepackage{gensymb}
\usepackage{epsfig}
\usepackage[left=90pt,right=90pt]{geometry}

\title{Comparison of surrogate-based uncertainty quantification methods for computationally expensive simulators}

\author{N. E. Owen\footnotemark[2]
\and P. Challenor\footnotemark[2]
\and P. P. Menon\footnotemark[2]
\and S. Bennani\footnotemark[3]}

\begin{document}
\maketitle

\newcommand{\slugmaster}{%
\slugger{juq}{xxxx}{xx}{x}{x--x}}
\renewcommand{\thefootnote}{\fnsymbol{footnote}}

\footnotetext[2]{College of Engineering, Mathematics and Physical Sciences, University of Exeter, Devon, UK.}
\footnotetext[3]{European Space Research and Technology Centre, TEC-ECN, Keplerlaan 1, 2201 AZ, Noordwijk, The Netherlands.}

\renewcommand{\thefootnote}{\arabic{footnote}}

\begin{abstract}
Polynomial chaos and Gaussian process emulation are methods for surrogate-based uncertainty quantification, and have been developed independently in their respective communities over the last 25 years. Despite tackling similar problems in the field, to our knowledge there has yet to be a critical comparison of the two approaches in the literature. We begin by providing a detailed description of polynomial chaos and Gaussian process approaches for building a surrogate model of a black-box function. The accuracy of each surrogate method is then tested and compared for two simulators used in industry: a land-surface model (\verb+adJULES+) and a launch vehicle controller (\verb+VEGACONTROL+). We analyse surrogates built on experimental designs of various size and type to investigate their performance in a range of modelling scenarios. Specifically, polynomial chaos and Gaussian process surrogates are built on Sobol sequence and tensor grid designs. Their accuracy is measured by their ability to estimate the mean, standard deviation, exceedance probabilities and probability density function of the simulator output, as well as a root mean square error metric, based on an independent validation design. We find that one method does not unanimously outperform the other, but advantages can be gained in some cases, such that the preferred method depends on the modelling goals of the practitioner. Our conclusions are likely to depend somewhat on the modelling choices for the surrogates as well as the design strategy. We hope that this work will spark future comparisons of the two methods in their more advanced formulations and for different sampling strategies.
\end{abstract}

\section{Introduction}
\label{Introduction}

Computer simulation of physical systems is now ubiquitous in science, because experimentation in the field can be expensive, time-consuming, or even impossible in practice. Examples include: regional ocean model systems\footnote{\texttt{http://myroms.org/}} \citep{Haidvogel:etal:2008}, modern global three/four-dimensional climate models \citep{HadGEM2:2011}, models for complex real world engineering problems (automobile, aerospace and construction) generated using general-purpose finite element programs\footnote{\texttt{http://www.lstc.com/applications}} \citep{Kirkpatrick:2000} and high fidelity mathematical models used in multi disciplinary design optimisation problems (for example wing design for a high speed civil transport aircraft) \citep{Giunta:etal:1999}. The computational models used for these purposes --- which we call \emph{simulators} --- are implemented as computer programs, and typically comprise of a system of (partial) differential equations involving a large number of input and output parameters. Due to the complexity of the equations which make up the model as well as the numerical tools needed to obtain their solution, it can take a substantial amount of time (hours, days, weeks) to execute a single run of the simulator. Running a simulator at various input parameter configurations to learn about a physical system is known as a computer experiment \citep{Sacks:etal:1989b}. \\

The field of uncertainty quantification in computer experiments has grown rapidly in recent years, where we aim to identify and reduce uncertainties found in all aspects of the computer modelling framework. Conventional approaches for tasks such as sensitivity analysis and calibration are based around Monte Carlo simulation.  However, exhaustively running an expensive simulator at a large number of input settings is simply not practical here. A solution to this problem is to instead use a smaller number of simulator runs to build a surrogate model, which is designed to provide an accurate and fast approximation to the simulator. In this way, the surrogate model can be used in place of the simulator for addressing uncertainty quantification objectives at a fraction of the original cost. Surrogate methods typically make no use of the equations that make up the simulator, assuming only a ``black-box'' model. Furthermore, they are usually applied to deterministic simulators which produce identical outputs if run at the same input settings. We also focus on simulators which are deterministic black-box functions in the present work. Two of the most popular and widely used surrogate methods, which are the focus of this work, are polynomial chaos and Gaussian process emulation. \\

In polynomial chaos, the simulator output is represented as a series expansion of functionals of the input parameters. Making use of orthogonality properties of polynomial families, the method uses information from the simulator to determine the coefficients of the expansion. Principally used by engineering and applied mathematics communities, the term ``polynomial chaos'' was first coined by \citet{Wiener:1938} who studied decompositions of Brownian motion. Polynomial chaos was initially applied to computer experiments in the seminal work of \citet{Ghanem:Spanos:1991b}, who were primarily interested in structural reliability problems. This early work in the field used chaos expansions made up of Hermite polynomials and Gaussian random variables. An extension of this was given by \citet{Xiu:Karniadakis:2003}, who incorporated non-Gaussian variables and polynomials from the Askey scheme in an approach known as \emph{generalized polynomial chaos}. Much of the early work implemented polynomial chaos for problems governed by known differential equations --- e.g., the Navier-Stokes equations in computational fluid dynamics \citep{LeMaitre:etal:2001,Knio:LeMaitre:2006} --- finding the expansion coefficients by solving an additional coupled system through a Galerkin projection. This so-called \emph{intrusive} approach requires knowledge of the equations that make up the simulator and thus cannot be applied for all computer experiments. Recently there has been a surge in the development of \emph{non-intrusive} alternatives, which evaluate the expansion coefficients using only repeated executions of the simulator. A big advantage of the non-intrusive class of methods is that the simulator is viewed as a black-box --- we can use already compiled code and do not need to tamper with the simulator itself (useful for legacy codes). Non-intrusive approaches are wide-ranging, and include sampling \citep{Ghiocel:Ghanem:2002,Reagan:etal:2003}, quadrature \citep{LeMaitre:etal:2002,Debusschere:etal:2004,Nobile:etal:2008a} and regression \citep{Berveiller:etal:2006,Blatman:Sudret:2011} based techniques. Another important example is the method of stochastic collocation \citep{Xiu:Hesthaven:2005} which uses interpolating polynomials. More advanced formulations of polynomial chaos have also been developed recently including: incorporating arbitrary probability distributions \citep{Wan:Karniadakis:2005}, multi-element \citep{Wan:Karniadakis:2006,Jakeman:etal:2013}, gradient-enhanced \citep{Liu:etal:2015}, adaptive and sparse \citep{Blatman:Sudret:2010a,Winokur:etal:2016}, oscillatory \citep{Witteveen:etal:2008}, multi-level and multi-index Monte Carlo hybrids \citep{Nobile:Tesei:2015,HajiAli:etal:2016} and multi-fidelity \citep{Narayan:etal:2014}. \\

By contrast, Gaussian process emulation treats the simulator as an unknown function of its inputs, which is modelled as a realisation of a stochastic process. Prior information about the simulator is expressed as a Gaussian process with mean and covariance functions, which are updated using runs of the simulator in an experimental design. Strongly related to Kriging methodology in geostatistics \citep{Cressie:2015}, Gaussian process emulation is primarily used by the statistics and applied mathematics communities and was first applied to the field of computer experiments by \citet{Sacks:etal:1989b}. Taking a frequentist viewpoint, their methodology used a regression model in conjunction with a zero-mean Gaussian process to best predict simulator output at an untried input configuration. The concept was soon viewed from the Bayesian perspective by \citet{Currin:etal:1991}, who were also concerned with prediction. Since this early work, the Bayesian methodology has been further refined \citep{Haylock:OHagan:1996,Kennedy:OHagan:2001} in terms of mathematical elegance, and to incorporate additional sources of uncertainty (for example model discrepancy, parametric uncertainty, observation error) found in computer experiments. Emulators have also been applied to address various objectives, including uncertainty analysis \citep{Oakley:OHagan:2002}, sensitivity analysis \citep{Oakley:OHagan:2004}, calibration \citep{Kennedy:OHagan:2001,Higdon:etal:2004} and history matching \citep{Williamson:etal:2013}. Similar to polynomial chaos, more advanced formulations of Gaussian process emulation have been developed including: multi-fidelity \citep{Kennedy:OHagan:2000,Forrester:etal:2007,LeGratiet:2013}, nested and hierarchical \citep{Oughton:Craig:2016}, sequential and adaptive \citep{Busby:2009,Loeppky:etal:2010}, gradient-enhanced \citep{Dwight:ZhongHua:2009} and dynamical or multivariate \citep{Conti:OHagan:2010,Fricker:etal:2013,Picheny:Ginsbourger:2013}. \\

It is clear that polynomial chaos and Gaussian process emulation are contrasting but related methods that both provide efficient surrogate-based uncertainty quantification for expensive simulators. Developed by separate communities but applied to computer experiments over roughly the same time period, to our knowledge there has yet to be a critical comparison of the two methods in the literature. The respective communities are clearly aware of each other; for instance \citet{OHagan:2013} gave a tutorial on polynomial chaos from a statistician's perspective. Furthermore, \citet{DiazDelaO:Adhikari:2011} proposed a method to reduce the computational cost of polynomial chaos by combining the expansion with an emulator, albeit in the field of stochastic finite elements. In more recent work, \citet{Liu:etal:2015} compare gradient-enhanced surrogate methods (including polynomial chaos and Gaussian process emulators) to a quasi-Monte Carlo approach for a geometry-induced uncertainty example in aerodynamics. \\

The purpose of this paper is to bring together the polynomial chaos and Gaussian process emulation communities --- in the spirit of \citet{OHagan:2013}, who recast polynomial chaos methodology in a statistical framework --- by providing a comprehensive description and comparison of the two techniques, assessing their relative advantages and disadvantages in a range of modelling scenarios. We concentrate on describing and comparing the approaches on a basic level, in the sense that although we are aware of more advanced techniques associated with both methods (for example multi-fidelity \citep{LeGratiet:2013,Narayan:etal:2014} or adaptive \citep{Blatman:Sudret:2010a,Loeppky:etal:2010} extensions), we believe a comparison of ``off-the-shelf'' methods would be more useful to practitioners in the area. Naturally, our particular modelling choices in constructing polynomial chaos and Gaussian process surrogates will have an effect on the conclusions we draw from our experiments. We hope that this work will spark further research into how the methods compare in their more advanced formulations. \\

In the examples we present, we seek to non-intrusively build polynomial chaos and Gaussian process surrogates to best represent a simulator output across the design space of the input parameters (the prediction objective outlined above), rather than tailor the surrogate for any specific uncertainty quantification task such as calibration or sensitivity analysis. This is to build a comparison of the two methods from the ground up, and we hope that the present paper will fuel further studies in other uncertainty quantification topics. We are primarily interested in the comparative accuracy of each method --- how well the surrogates approximate the simulator --- with changes in design type and size. We achieve this by building and validating polynomial chaos and Gaussian process surrogates on land-surface (\verb+adJULES+) and launch control vehicle (\verb+VEGACONTROL+) simulators. For all cases, the surrogates are built and compared on exactly the same design size and type for a fair analysis. Even so, we acknowledge that our sampling strategy is likely to have some effect on the performance of the two surrogate methods. We hope that future research will be able to investigate the impact of design choice on the comparison. \\

The outline of the remainder of the paper is as follows. In \S\ref{Surrogates} we introduce notation for surrogate modelling and give detailed descriptions of polynomial chaos (\S\ref{PC}) and Gaussian process emulation (\S\ref{GP}) methodologies. We also outline methods for validation and comparison of surrogate models in \S\ref{Validation}. In \S\ref{Application} we give details of our experiments and present results from the \verb+adJULES+ (\S\ref{adJULES}) and \verb+VEGACONTROL+ (\S\ref{VEGACONTROL}) simulators. We conclude the paper with some discussion in \S\ref{Discussion}. 

\section{Surrogate modelling}
\label{Surrogates}

We proceed by introducing our notation for surrogate modelling. Consider a deterministic black-box simulator, $\eta(\cdot)$, as a function of $n$ input parameters, $\mathbf{x} = \left\{x_1,\ldots,x_n \right\} \subset \mathbb{R}^n$, producing a set of outputs $\mathbf{y}$. Suppose we are concerned with only one of those outputs, $y \in \mathbb{R}$. As outlined in \S\ref{Introduction}, we assume that the simulator is computationally expensive, in that the evaluation of simulator output $y^{(i)}$ from a single input configuration $\mathbf{x}^{(i)} = \{x_1^{(i)},\ldots,x_n^{(i)}\}$ may take a substantial amount of time. Computational and time issues mean that we are restricted to making a small number of runs of the simulator, which we use to build a surrogate model. \\

Let $\mathcal{D} = \{\mathbf{x}^{(1)},\ldots,\mathbf{x}^{(m)}\}$ be a set of input configurations in an experimental design, and $\mathcal{Y} = \{y^{(1)},\ldots,y^{(m)}\}^{T}$ be the corresponding output obtained by evaluating $y^{(i)} = \eta(\mathbf{x}^{(i)}), \, i = 1,\ldots,m$. The experimental design should follow some basic principles \citep{Santner:etal:2003} to ensure an accurate representation of the simulator. For surrogate modelling the experimental design should generally be space-filling in the input region of interest, and Latin Hypercube \citep{McKay:etal:1979} or tensor (factorial) grid designs are popular here. Of course, there is no need to run a deterministic simulator at the same input configuration more than once, so design points should be strategically placed to maximise the information gained from the number of simulator runs you can afford. The design may also be restricted because it is itself an integral part of building the surrogate --- for example, if a quadrature rule is required in the fitting process. \\

The surrogate model, denoted $\hat{\eta}(\cdot)$, is built non-intrusively using only information contained in $\mathcal{D}$ and $\mathcal{Y}$. Once built, it should provide a fast approximation to the simulator at any untried input configuration $\mathbf{x}^{(*)}$ so that it may be used in place of it for subsequent uncertainty quantification tasks. Other desirable (but not necessary) properties for a surrogate are for it to \emph{exactly} reproduce simulator output at the experimental design points (since the simulator is deterministic), and for it to provide some form of uncertainty measure along with its predictions.

\subsection{Polynomial chaos}
\label{PC}

In polynomial chaos, the simulator output is modelled as a function of the input parameters in a series expansion. Since we are uncertain about the correct configuration of the input parameters, in the polynomial chaos framework we begin by representing them as a random vector $\mathbf{X}$ with joint probability density function (PDF) $f_{\mathbf{X}}(\mathbf{x})$. We are primarily interested in the induced uncertainty on the simulator output, $Y \equiv \eta(\mathbf{X})$. In the following, we assume independence of the input parameters such that their joint PDF may be written as
$f_{\mathbf{X}}(\mathbf{x}) = \prod_{j=1}^n f_{X_j}(x_j)$, where $f_{X_j}(x_j)$ is the marginal PDF of input parameter $X_j$. For computer experiments, distributions for the input parameters can usually be found or estimated through expert elicitation. In the case of parameter dependence or a random field, an extra decorrelation step is required (for example, using a Rosenblatt transformation or a Karhunen-Lo\`{e}ve expansion) \citep{Eldred:etal:2008}.  \\

We assume that the uncertain simulator output $Y$ is a second-order stationary process; that is, it can be characterised by its first and second moments. It can be shown \citep{Xiu:Karniadakis:2003} that $Y$ can expanded onto an orthogonal polynomial basis as follows, 

\begin{equation}
\label{PCE}
Y \equiv \eta(\mathbf{X}) = \sum_{\boldsymbol{\alpha} \in \mathbb{N}^n} a_{\boldsymbol{\alpha}} \boldsymbol{\psi}_{\boldsymbol{\alpha}}(\mathbf{X}) \, .
\end{equation}

The series in \eqref{PCE} is known as a \emph{polynomial chaos expansion}. The $a_{\boldsymbol{\alpha}}$'s are unknown expansion coefficients (to be determined) and the $\boldsymbol{\psi}_{\boldsymbol{\alpha}}(\mathbf{X})$'s are known multivariate polynomials (specified a priori). The multidimensional summation index $\boldsymbol{\alpha}$ will be described in more detail later. The selection of multivariate polynomials used in \eqref{PCE} depends on a number of orthogonality properties, which we will now outline. \\

Consider the set $\{\psi_k^{(j)}, k \in \mathbb{N}\}$ a family of polynomials in terms of the parameter $X_j$, where $k$ denotes the polynomial degree. For example, $\psi_5^{(2)}$ denotes a fifth order polynomial in terms of the parameter $X_2$. A set of polynomials are said to be \emph{orthogonal} with respect to a probability distribution $f_{X_j}$ if the following inner product holds \citep{Ghanem:Spanos:1991b}:

\begin{equation}
\label{orthogonality1}
\left\langle \psi_{k}^{(j)}(X_j), \psi_{l}^{(j)}(X_j) \right\rangle \equiv \mathbb{E}\left[\psi_{k}^{(j)}(X_j) \, \psi_{l}^{(j)}(X_j) \right] \equiv \int_{\mathcal{X}_j} \psi_{k}^{(j)}(x_j) \, \psi_{l}^{(j)}(x_j)\,  f_{X_j}(x_j) \, \mathrm{d}x = \gamma_k^2 \, \delta_{kl} \, ,
\end{equation} 

where $\mathcal{X}_j$ is the support of $f_{X_j}$ and $\delta_{kl}$ is the Kronecker delta. The normalisation constants $\gamma_k$ are unique to the chosen family of polynomials and are known in practice. \\

Since we have assumed independence of the input parameters, the multivariate polynomials used in \eqref{PCE} are constructed using a tensor product of the $n$ families of univariate polynomials,

\begin{equation*}
\boldsymbol{\psi}_{\boldsymbol{\alpha}}(\mathbf{x}) = \psi_{\alpha_1}^{(1)}(x_1) \times \cdots \times \psi_{\alpha_n}^{(n)}(x_n) \, ,
\end{equation*}

where we use the multidimensional index notation $\boldsymbol{\alpha} \equiv \{\alpha_1,\ldots,\alpha_n\} \in \mathbb{N}^n$ to define the degree of each one-dimensional polynomial. As an example, $\boldsymbol{\psi}_{\{2,4,3\}}$ would correspond to the multivariate polynomial $\psi_2^{(1)}\psi_4^{(2)}\psi_3^{(3)}$. The multivariate version of the inner product in \eqref{orthogonality1} is

\begin{equation}
\label{orthogonality2}
\left\langle \psi_{\boldsymbol{\alpha}}(\mathbf{X}), \psi_{\boldsymbol{\beta}}(\mathbf{X}) \right\rangle \equiv \mathbb{E}\left[\psi_{\boldsymbol{\alpha}}(\mathbf{X}) \, \psi_{\boldsymbol{\beta}}(\mathbf{X}) \right] \equiv \int_{\mathcal{X}} \psi_{\boldsymbol{\alpha}}(\mathbf{X}) \, \psi_{\boldsymbol{\beta}}(\mathbf{X})\,  f_{\mathbf{X}}(\mathbf{x}) \, \mathrm{d}\mathbf{x} = \gamma_{\boldsymbol{\alpha}}^2 \, \delta_{\boldsymbol{\alpha}\boldsymbol{\beta}} \, ,
\end{equation}

where $\mathcal{X} = \mathcal{X}_1 \times \cdots \times \mathcal{X}_n$ is the support of $f_{\mathbf{X}}(\mathbf{x})$. \\

As a consequence of \eqref{orthogonality1} and \eqref{orthogonality2}, the polynomial terms used in the expansion are entirely dependent on the probability distributions of the input parameters. As mentioned in \S\ref{Introduction}, early work involved the use of Gaussian random variables for the input parameters, which results in Hermite polynomials being used in the expansion \citep{Wiener:1938,Ghanem:Spanos:1991b}. In the extended approach of generalized polynomial chaos \citep{Xiu:Karniadakis:2003}, most common types of  distributions for the inputs are accounted for, resulting in polynomials from the Askey scheme being used in the polynomial chaos expansion. To give some examples, if the input parameters have a Gamma distribution then we use Laguerre polynomials, whereas Poisson distributed inputs are coupled with Charlier polynomials \citep{Xiu:Karniadakis:2003}. A mix of distributions across the input parameters is also supported, and in this case the expansion would be made up of a mix of the relevant polynomial families. For the examples used in this paper, we assume all input parameters are uniformly distributed which gives rise to the use of Legendre polynomials. This is for two reasons: firstly we have no information to say they are distributed otherwise, and secondly it leads to a fairer comparison with Gaussian process methods (since we give no advantage to polynomial chaos by hand-picking the most appropriate polynomial families). \\

To be used in practice, the polynomial chaos expansion in \eqref{PCE} must be truncated after $N$ terms. It is common to keep polynomials $\boldsymbol{\psi}_{\boldsymbol{\alpha}}$ with degree up to $p$:

\begin{equation}
\label{truncPCE}
Y \approx \eta_p(\mathbf{X}) = \sum_{0 \leq |\boldsymbol{\alpha}| \leq p} a_{\boldsymbol{\alpha}} \boldsymbol{\psi}_{\boldsymbol{\alpha}}(\mathbf{X}) \equiv \mathbf{a}^{T} \Psi(\mathbf{X}) \, ,
\end{equation}

where $\mathbf{a}$ and $\Psi(\mathbf{X})$ are vectors containing the coefficients and polynomial elements of the truncated expansion respectively. The parameter $p$ is usually referred to as the truncation order. The truncation depends on how you define $|\boldsymbol{\alpha}|$, and two truncation schemes are used here:

\begin{enumerate}
\item $|\boldsymbol{\alpha}| = \sum_{j=1}^n \alpha_j$
\item $|\boldsymbol{\alpha}| = $ each $\alpha_j, j=1,\ldots,n$, in $\boldsymbol{\alpha}$ considered individually 
\end{enumerate}

In this paper, (i) is referred to as a \emph{total order} truncation since we include all polynomials where summation of the multidimensional index is less than or equal to the truncation order $p$. In contrast, (ii) is referred to as a \emph{tensor product} truncation, since the restriction of order is applied in each dimension and all combinations of possible one-dimensional polynomials are included. To give a worked example, suppose we require a quadratic expansion ($p=2$) in terms of two input parameters ($n=2$). A total order expansion would include the polynomial  $\Psi_{\{0,0\}} = \psi_0^{(1)}\psi_0^{(2)}$ as well as $\psi_1^{(1)}\psi_0^{(2)}$, $\psi_0^{(1)}\psi_1^{(2)}$, $\psi_2^{(1)}\psi_0^{(2)}$, $\psi_1^{(1)}\psi_1^{(2)}$ and $\psi_0^{(1)}\psi_2^{(2)}$, since the multidimensional index $\boldsymbol{\alpha} = \{\alpha_1,\alpha_2\}$ sums to $p=2$ or less in each case. The alternative tensor product expansion would also include the polynomials $\psi_2^{(1)}\psi_1^{(2)}$, $\psi_1^{(1)}\psi_2^{(2)}$ and $\psi_2^{(1)}\psi_2^{(2)}$ since the restriction of order is applied in each input dimension separately. \\

The number of terms, $N$, for the truncated polynomial chaos expansion in \eqref{truncPCE} grows rapidly with the number of input parameters $n$ and the truncation order $p$. Specifically, $N = {n +p \choose p}$ for a total order expansion and $N = (p+1)^n$ for a tensor product expansion. Since $n$ is typically fixed for the simulator, one must choose $p$ in accordance with computational restraints, that is, the number of simulator runs you can afford. To successfully fit a polynomial chaos expansion, at least $N$ design points are required to avoid an under-determined system. An alternative truncation scheme has recently been proposed in the method of sparse polynomial chaos \citep{Blatman:Sudret:2010a,Blatman:Sudret:2011}. Here, the curse of dimensionality is reduced by using an algorithm to select the most influential terms in the expansion. However, for simplicity this will not be considered here. \\

% Possible figure for N_t
%\begin{figure}[htb]
%\centering 
%\includegraphics[width=0.8\textwidth]{./Plots/PCplot.png}
%\caption{ghj.}
%\end{figure}

With the initial selection of polynomial terms complete, all that remains is to estimate the expansion coefficients $a_{\boldsymbol{\alpha}}$ using information from the simulator. Methods to compute the expansion coefficients are classed intrusive and non-intrusive depending on whether they require knowledge of the equations that make up the simulator. Since we are treating the simulator as a black-box, we focus only on non-intrusive approaches. Specifically, the following two non-intrusive methods are used:

\begin{enumerate}
\item \emph{regression}: uses the least-squares solution to a linear system
\item \emph{spectral projection}: uses the numerical solution of an integral, e.g., using quadrature 
\end{enumerate} 

Other non-intrusive methods such as sampling and stochastic collocation are not considered here. \\

We now proceed with a description of the two non-intrusive methods. For the regression approach, the expansion coefficients are found by minimising \citep{Blatman:Sudret:2011}

\begin{equation*}
\hat{\mathbf{a}} = \mathrm{arg} \min_{\mathbf{a} \in \mathbb{R}^{N}} \left[ \sum_{i=1}^m \left(\eta(\mathbf{x}^{(i)}) -  \mathbf{a}^{T} \Psi(\mathbf{x}^{(i)}) \right)^2 \right] \, ,
\end{equation*}

that is, the sum of the squared difference between the simulator output and the truncated PCE in \eqref{truncPCE} at the experimental design points. The least-squares solution is

\begin{equation*}
\hat{\mathbf{a}} = \left(\boldsymbol{\Psi}^{T}\boldsymbol{\Psi}\right)^{-1}\boldsymbol{\Psi}^{T}\mathcal{Y} \, ,
\end{equation*}

where 

\begin{equation*}
\boldsymbol{\Psi}_{ij} = \boldsymbol{\psi}_{\boldsymbol{\alpha}_j}(\mathbf{x}^{(i)}), \qquad i = 1,\ldots,m, \, j=1,\ldots,N \, ,
\end{equation*}

is a data matrix containing the $N$ polynomial terms in \eqref{truncPCE}, evaluated at the $m$ design points. This approach is also known as \emph{point-collocation}, and the design size $m$ must be greater than or equal to the number of expansion terms $N$ to produce a well-conditioned linear system. \citet{Hosder:etal:2007} suggest $m = 2N$ (\emph{twice over-determined}), but we also experiment with $m = N$ (\emph{uniquely-determined}) in our experiments. This is to test the capability of the regression approach when faced with small design sizes. We acknowledge that for the uniquely-determined regression, or cases where $m<N$, sparse approximation techniques based on $l_1$-minimisation, LASSO or least-angle regression (for example, \citet{Blatman:Sudret:2010a,Blatman:Sudret:2011}) are likely to give more stable results. However, since we seek to perform a comparison off ``off-the-shelf'' methods, we only consider the more common least-squares method. Finally, this approach is combined with a polynomial chaos expansion with total order truncation as is common in the literature \citep{Eldred:Burkardt:2009}. \\

The spectral projection approach projects the simulator output against each polynomial term using inner products and employs orthogonality to extract each expansion coefficient. Taking the inner product of \eqref{truncPCE} with respect to $\boldsymbol{\psi}_{\boldsymbol{\alpha}}$ and enforcing orthogonality \citep{Eldred:etal:2008} yields the following expression for the coefficients,

\begin{equation}
\label{projection}
a_{\boldsymbol{\alpha}} = \frac{1}{\gamma_{\boldsymbol{\alpha}}^2} \langle Y, \psi_{\boldsymbol{\alpha}}(\mathbf{X})\rangle = \frac{1}{\gamma_{\boldsymbol{\alpha}}^2} \int_{\mathcal{X}} Y \, \psi_{\boldsymbol{\alpha}}(\mathbf{X}) \, f_{\mathbf{X}}(\mathbf{x}) \, \mathrm{d}\mathbf{x} \, .
\end{equation}

We cannot evaluate this expression analytically since we do not know the simulator output at all points in the design space. To estimate the expansion coefficients we need to use a numerical method to approximate the integral on the right hand side of \eqref{projection}. Early work used Monte Carlo or sampling based methods \citep{Ghiocel:Ghanem:2002,Reagan:etal:2003}, although improvements can be made through the use of quadrature methods for numerical integration \citep{LeMaitre:etal:2002,Debusschere:etal:2004,Eldred:etal:2008}. In this case, the experimental design $\mathcal{D}$ should encompass a quadrature rule to ensure accurate approximation of the integral. There are several choices of quadrature abscissae, but in our experiments we implement a tensor product of Gauss-Legendre points to complement the use of multivariate Legendre polynomials in \eqref{truncPCE} \citep{Eldred:etal:2008}. Gaussian-type quadrature rules are the natural choice for our ``off-the-shelf'' comparison as they are standard for polynomial chaos, accounting for the probability weight in the integral in \eqref{projection}. However, other choices of quadrature rules (for example, Clenshaw-Curtis) may perform better in certain scenarios. A tensor product truncation scheme for the polynomial chaos expansion is also used to align with the tensor product experimental design. A disadvantage of tensor product quadrature rules is that they suffer from the curse of dimensionality, because the number of design points required grows rapidly with the input dimension. When the input dimension becomes large (e.g., $n > 5$), sparse grid quadrature (not to be confused with sparse polynomial chaos) is an alternative method which can reduce the computational burden while retaining high accuracy \citep{Xiu:Hesthaven:2005,Nobile:etal:2008a}. However, the input dimension in our examples is relatively small so sparse grids are not required. \\

Finally the estimated expansion coefficients, $\hat{a}_{\boldsymbol{\alpha}}$, are substituted into the truncated polynomial chaos expansion in \eqref{truncPCE}. The resulting polynomial chaos surrogate can be used to approximate simulator output at any untried input configuration. In this way, features of the simulator output, such as its mean and variance, may be estimated by using the surrogate in a Monte Carlo fashion.

\subsection{Gaussian process emulation}
\label{GP}

In the Gaussian process emulation framework, the simulator is viewed as an unknown function which we model as a realisation of a stochastic process. In particular, the main assumption we make is that the simulator can be represented by a Gaussian process; that is, simulator output obtained at different input configurations can be modelled with a joint Gaussian distribution. \\

A Gaussian process is fully specified by a mean function, $M(\cdot)$, and a covariance function, $V(\cdot,\cdot)$. Taking a Bayesian perspective, we put priors on these to reflect any knowledge or information we may have about the simulator. As is common in the literature, we restrict ourselves to mean functions that can be represented as a sum of basis functions, $M(\mathbf{x}) = h(\mathbf{x})^T \boldsymbol{\beta}$. Here, $h(\mathbf{x})$ is a $q \times 1$ vector of known basis or regression functions and $\boldsymbol{\beta}$ are the corresponding coefficients, to be estimated. The mean function should describe any known global behaviour of the simulator output across the input space, and it is common to use a simple constant or linear trend in the absence of any information. In our experiments, however, we determine the regression functions using a preliminary regression analysis in conjunction with a stepwise algorithm. In this way, we allow data from the simulator runs to determine a suitable form for the mean function. \\

The covariance function describes how we expect the simulator output to be correlated as a function of two input configurations, say $\mathbf{x}$ and $\mathbf{x}'$. Where the mean function models the global behaviour of the simulator output across the input space, the covariance function controls local behaviour. Usual practice is to define $V(\mathbf{x},\mathbf{x}') = \lambda^2 C(\mathbf{x},\mathbf{x}';\boldsymbol{\delta})$, where $\lambda^2$ is the Gaussian process variance, and $C(\mathbf{x},\mathbf{x}';\boldsymbol{\delta})$ is a correlation function dependent on parameters $\boldsymbol{\delta} = \{\delta_1,\ldots,\delta_n\}$. These parameters are referred to as \emph{correlation lengths}, and control the strength of the correlation in  output for each of the $n$ input parameters. In general, a smaller correlation length corresponds to a rougher process and a larger correlation length corresponds to a smoother process. The Gaussian process variance describes the extent to which the simulator can deviate away from the mean function. \\

We focus on continuous and stationary correlation functions, where the correlation is simply a function of the Euclidean distance between two inputs. We define this distance between $\mathbf{x}$ and $\mathbf{x}'$ separately along each input dimension $j$ as $d_j \equiv |x_j - x_j'|$, $j = 1,\ldots,n$. The most widely used correlation function in the literature, which we use in our examples, is the squared exponential (sometimes referred to as Gaussian) \citep{Rasmussen:Williams:2006}. The squared exponential covariance function has the following form for $n$ input parameters,

\begin{equation*}
C(\mathbf{x},\mathbf{x}') = \exp \left( -\frac{1}{2} \sum_{j=1}^n \left(\frac{d_j}{\delta_j}\right)^2\right) \, .
\end{equation*} 

This correlation function is infinitely differentiable, which results in a Gaussian process that is very smooth. Since this may not necessarily be the case for the simulators in our examples, we also experiment with a Mat\'{e}rn covariance function \citep{Rasmussen:Williams:2006}. This class of covariance function has an extra scale parameter, $\nu$, which controls the differentiability of the process. Specifically, the process is $\left\lfloor{\nu}\right\rfloor$ times differentiable \citep{Rasmussen:Williams:2006}. Commonly, this parameter is set to a half-integer to simplify the form of the covariance function. We use $\nu=5/2$, which results in the the following form for the Mat\'{e}rn covariance function in $n$ input parameters,

\begin{equation*}
C(\mathbf{x},\mathbf{x}') = \left[\prod_{j=1}^n \left(1+\frac{\sqrt{5}d_j}{\delta_j}+\frac{5 d_j^2}{3 \delta_j^2} \right) \right] \exp \left(-\sqrt{5} \sum_{j=1}^n \frac{d_j}{\delta_j} \right) \, .
\end{equation*}

This particular covariance function results in a Gaussian process that is twice differentiable, and consequently a rougher process than if using a squared exponential covariance function. \\

The complete Gaussian process prior for the simulator can be written

\begin{equation*}
\eta(\mathbf{x}) \, | \, \boldsymbol{\beta},\lambda^2,\boldsymbol{\delta} \sim \mathcal{GP}(h(\mathbf{x})^T \boldsymbol{\beta},\lambda^2 C(\mathbf{x},\mathbf{x};\boldsymbol{\delta})) \, . 
\end{equation*}

This prior process is very flexible and can represent a wide range of simulator output behaviours. As with our construction of the polynomial chaos expansion, the choices we have made in specifying the Gaussian process are standard in the literature. More advanced formulations (for example, nonstationary covariance functions) are likely to perform better in certain scenarios, but would not make for a balanced comparison. \\

The unknown parameters $\boldsymbol{\beta}$, $\lambda^2$ and $\boldsymbol{\delta}$ must be estimated using information from the simulator runs, $\mathcal{Y}$. If we assume a standard non-informative prior, namely $p(\boldsymbol{\beta},\lambda^2) \propto \lambda^{-2}$, an analytical marginalisation is possible of the $\boldsymbol{\beta}$ and $\lambda^2$ parameters \citep{Haylock:OHagan:1996}. This results in the following marginal likelihood for the correlation lengths

\begin{equation}
\label{likelihood}
L(\boldsymbol{\delta}|\mathcal{Y}) \propto (\hat{\lambda}^2)^{-(m-q)/2}|A|^{-1/2} |H^{T} A^{-1} H|^{-1/2} \, ,
\end{equation}

where $\hat{\lambda}^{2} = (\mathcal{Y}-H\hat{\boldsymbol{\beta}})^{T}A^{-1}(\mathcal{Y}-H\hat{\boldsymbol{\beta}})$ and $\hat{\boldsymbol{\beta}} = (H^{T}A^{-1}H)^{-1}H^{T}A^{-1}\mathcal{Y}$. Furthermore, the experimental design points define the regression matrix $H$, with $i^{\mathrm{th}}$ row $H_i = h(\mathbf{x}^{(i)})^{T}$, and the correlation matrix $A$, which has elements $A_{ij} = c(\mathbf{x}^{(i)},\mathbf{x}^{(j)})$. The use of a non-informative prior for the regression and variance parameters means that some accuracy is lost in their estimation, but in practice this is not noticeable. It also leads to a fairer comparison to polynomial chaos methods since we are not giving Gaussian process emulation any advantage. \\

The correlation length parameters can be estimated in a number of ways. A fully Bayesian analysis would proceed with a numerical marginalisation of $\boldsymbol{\delta}$, e.g., using a MCMC algorithm. This is typically a computationally expensive approach and is rarely used in practice. We favour the method outlined in \citet{Kennedy:OHagan:2001}, who simply find a maximum likelihood estimate of $\boldsymbol{\delta}$ and assume it as the true value. Of course, this `plug-in' approach does not fully account for the uncertainty in these parameters, but the loss of accuracy is minimal and the savings in computational cost typically outweigh this disadvantage \citep{Oakley:OHagan:2004}. We estimate the correlation lengths by maximising the logarithm of the likelihood in \eqref{likelihood},

\begin{equation}
\label{loglikelihood}
\hat{\boldsymbol{\delta}} = \mathrm{arg} \max_{\boldsymbol{\delta}} \left[\log L(\boldsymbol{\delta}|\mathcal{Y}) \right] \, .
\end{equation} 

In our experiments we use the \verb+R+ package \verb+DiceKriging+ \citep{DiceKriging} to fit Gaussian process emulators. This maximises the log-likelihood in \eqref{loglikelihood} using the ``L-BFGS-B'' optimisation algorithm, which is a quasi-Newton scheme constrained by given lower and upper bounds \citep{Byrd:etal:1995,Park:Baek:2001}. This algorithm makes use of analytical gradients of the log-likelihood with respect to the correlation length parameters to speed up the optimisation. The number of iterations required for accurate approximation of the correlation lengths is problem dependent, but note that the main cost stems from inverting the $m \times m$ correlation matrix $A$. Strategies for using Gaussian process emulators for large datasets are presented in \citet{Rasmussen:Williams:2006}. \\  

After the correlation lengths have been estimated the result is a Gaussian process surrogate model, known as an \emph{emulator}. The emulator specifies a full posterior distribution for its predictions which is a Student's-$t$ distribution with $m-q$ degrees of freedom \citep{Haylock:OHagan:1996,Kennedy:OHagan:2001},

\begin{equation}
\label{posterior}
\eta(\mathbf{x}) \, | \, \mathcal{Y},\boldsymbol{\delta} \sim t_{m-q} (M^{*}(\mathbf{x}),V^{*}(\mathbf{x},\mathbf{x})) \, .
\end{equation}

Here we have introduced the posterior mean function,

\begin{equation*}
M^{*}(\mathbf{x}) = h(\mathbf{x})^{T}\hat{\boldsymbol{\beta}} + R(\mathbf{x})A^{-1}(\mathcal{Y}-H\hat{\boldsymbol{\beta}}) \, ,
\end{equation*}

and posterior covariance function, 

\begin{equation*}
V^{*}(\mathbf{x},\mathbf{x}') = \frac{\hat{\lambda}^2}{m-q-2} \left[c(\mathbf{x},\mathbf{x}') - R(\mathbf{x}) A^{-1}	R(\mathbf{x}')^{T} + Q(\mathbf{x}) (H^{T}A^{-1}H)^{-1}  Q(\mathbf{x}')^{T} \right] \, ,
\end{equation*}

where $Q(\mathbf{x}) = h(\mathbf{x})^{T} - R(\mathbf{x})A^{-1}H$, and $R(\mathbf{x})$ is a vector of correlations between $\mathbf{x}$ and the experimental design $\mathcal{D}$. A key feature of the Gaussian process emulator is that it interpolates the simulator output at the experimental design points (without a nugget parameter, see \citet{Andrianakis:Challenor:2012}). At an untried input configuration, the emulator gives a prediction of simulator output by taking the maximum of the posterior. In contrast to polynomial chaos methods, the distribution in \eqref{posterior} means the emulator also readily provides uncertainty information. Uncertainty in the output at the experimental design points is zero, but grows as we move away from them. This can be helpful in identifying areas in the design space where the emulator is particularly confident, or areas where we may wish to make more simulator runs because uncertainty is high. Nevertheless, it is common to simply use emulator predictions as a fast approximation to the simulator, enabling surrogate-based uncertainty quantification in a Monte Carlo fashion.

\subsection{Validation and comparison of surrogate models}
\label{Validation}

Once we have built a surrogate model, it is important to assess its quality. We may check the surrogate's ability to represent the simulator using a number of validation exercises, which tend to fall into two categories. Firstly, we can make $m'$ additional runs of the simulator in a \emph{validation design}, say $\mathcal{D}' = \{\mathbf{x}^{(1')},\ldots,\mathbf{x}^{(m')}\}$. The simulator output at the validation design, $\eta(\mathbf{x}^{(i)})$, $i = 1',\ldots,m'$, can then be compared to corresponding surrogate predictions at these points using various metrics. It is important that these design points are distinct from the original experimental design $\mathcal{D}$ to fully test the quality of the surrogate. Of course, this approach requires the extra computational cost of obtaining simulator output at the validation design. If this is not possible, a second approach is to re-use the original experimental design points in some way. The main class of methods here are based on cross-validation exercises, and leave-one-out cross validation is a popular approach. This is where surrogates are sequentially rebuilt on the original experimental design with one point withheld, and then used to predict simulator output at the withheld input configuration. A benefit of this approach is that no extra simulator runs are required, but computational cost does arise from refitting the surrogate model for the $m$ sub-designs. In this paper we implement the first of the outlined validation approaches, since for the simulators used in our examples it is relatively cheap to produce a set of validation runs. Furthermore, cross-validation is not compatible with the quadrature approach for polynomial chaos since we must use all quadrature abscissae for the method to work. \\

Given simulator output at the validation design and corresponding surrogate predictions, a  number of metrics can be used to assess performance. In the context of analysing two contrasting surrogate approaches, it is vital that we use metrics that compare Gaussian process emulation and polynomial chaos methods in an unbiased and fair way. Validation metrics specific to either approach, for example those for Gaussian process emulators outlined in \citet{Bastos:OHagan:2009} which take correlated predictions into account, are not necessarily suited to such a comparison. We proceed by presenting the validation metrics used in our experiments, which are computed using a $1000$-point Latin Hypercube \citep{McKay:etal:1979} validation design, independent to the experimental designs used to build the surrogates. Firstly, we use root mean square error (RMSE),

\begin{equation}
\label{rmse}
RMSE = \sqrt{\frac{1}{m'}\sum_{i=1'}^{m'} \left(\eta(\mathbf{x}^{(i)})-\hat{\eta}(\mathbf{x}^{(i)})\right)^2} \, ,
\end{equation} 

to assess the accuracy of the surrogate across the input space. Since this in effect measures the distance between the surrogate and the simulator, a lower RMSE is favourable. Secondly, let $\mu_Y$ and $\sigma_Y$ denote the true values of the simulator output mean and variance respectively. We estimate these using the surrogate in the following empirical fashion,

\begin{equation}
\label{mean}
\mu_Y \approx \hat{\mu}_Y = \frac{1}{m'}\sum_{i=1'}^{m'} \hat{\eta}(\mathbf{x}^{(i)}) \, ,
\end{equation}

\begin{equation}
\label{sd}
\sigma_Y \approx \hat{\sigma}_Y = \sqrt{\frac{1}{m'}\sum_{i=1'}^{m'}(\hat{\eta}(\mathbf{x}^{(i)})-\hat{\mu}_Y)^2} \, .
\end{equation}

The mean and standard deviations estimated by the surrogates can be compared to simulator quantities also computed at the validation design points. We acknowledge that for polynomial chaos, estimates of these quantities can be derived analytically from the expansion \citep{Xiu:Karniadakis:2003} but instead use the above Monte Carlo approach for fairer comparison to Gaussian process emulation. Thirdly, we are interested in the capability of the surrogates to estimate the chance of rare events, represented by the exceedance probabilities $\mathrm{Pr}(Y \geq \mu_Y + \kappa \sigma_Y)$. We estimate this quantity using the surrogate as follows,

\begin{equation}
\label{exceed}
\mathrm{Pr}(Y \geq \mu_Y + \kappa \sigma_Y) \approx \frac{1}{m'}\sum_{i=1'}^{m'} \boldsymbol{1}\left(\hat{\eta}(\mathbf{x}^{(i)}) \geq \mu_Y + \kappa \sigma_Y\right) \, ,
\end{equation}

where $\boldsymbol{1}$ denotes the indicator function, and we set $\kappa = 2, 3$ for probabilities of exceeding two and three standard deviations above the mean. Since the simulator has only been evaluated $1000$ times at the validation design, naturally we can only have approximate values of $\mu_Y$, $\sigma_Y$ and $\mathrm{Pr}(Y \geq \mu_Y + \kappa \sigma_Y)$. To account for the uncertainty on these quantities induced by finite simulator evaluations, we perform a bootstrap analysis to obtain $95\%$ bootstrap percentile intervals \citep{Efron:Tibshirani:1994}. Finally, we wish to test the surrogate's ability to reconstruct the PDF of the simulator output. We do this by smoothing the surrogate predictions at the validation design with a kernel density estimator. For all examples we use a Gaussian kernel and estimate the kernel bandwidth for the simulator output at the validation design using the method outlined in \citet{Silverman:1986}. The estimated bandwidth for the simulator is then kept the same between surrogate methods for consistency. \\

Since emulators provide uncertainty information, we also use the re-sampling approach outlined in \newline \citet{Oakley:OHagan:2002} to obtain $95\%$ confidence intervals about Gaussian process point estimates of all the validation metrics outlined above. This method works by sampling from the posterior distribution given in \eqref{posterior} a large number of times and averaging across the validation metrics computed for each sample. 

\section{Experiments}
\label{Application}

In this section we present results comparing polynomial chaos and Gaussian process emulation for two simulators, named \verb+adJULES+ and \verb+VEGACONTROL+. Our comparison uses the validation metrics outlined in \S\ref{Validation} to assess surrogate accuracy for different choices of experimental design. We investigate the performance of the two surrogate methodologies under changes in the type and size of experimental design used to build the surrogate. We proceed by describing the rationale behind our experimental design choices. A summary of the designs used in the \verb+adJULES+ and \verb+VEGACONTROL+ experiments is given in Table \ref{designs}. \\

We have grouped our experimental designs into three classes to distinguish between the different methods used to compute the coefficients of the polynomial chaos expansions. Recall in \S\ref{PC}, that two strategies were outlined for estimating the polynomial chaos expansion coefficients: \emph{regression} and \emph{spectral projection} (quadrature). For the regression approach there are no restrictions on the experimental design, but one should generally use a space-filling design with enough design points to ensure a well-determined system. For this purpose we use Sobol sequence designs \citep{Niederreiter:1988}, and these make up \emph{class 1} and \emph{class 2} designs (the difference between the two will be explained shortly). Sobol sequences are used (instead of the more popular Latin Hypercube designs) because of the fact that we can take subsets of the original design and it remains a Sobol sequence. In this way, we can save computational time by evaluating the simulator at one large Sobol sequence design and use its subsets to build surrogates based on smaller design sizes. For the spectral projection approach, numerical quadrature is used to evaluate the coefficients, therefore the experimental designs in this case are quadrature rules. As explained in \S\ref{PC}, we implement a tensor (factorial) grid of one-dimensional Gauss-Legendre quadrature rules, and these make up \emph{class 3} designs. \\

Each design class comprises four distinct designs of increasing size. Like the class of design, the design sizes are chosen according to how we construct the polynomial chaos expansion.  Recall that the design size ($m$) for polynomial chaos must be at least as big as the number of terms in the expansion ($N$), which depends on truncation order ($p$), the input parameter dimension ($n$) and the truncation scheme used (see \S\ref{PC} for the exact relation). Since $n$ is fixed for the simulators in our examples, we experiment with truncation orders of $p=1,2,3,4$ in our polynomial chaos expansions (we do not experiment with $p>4$ due to computational restraints). This fixes the number of terms in each of the polynomial chaos expansions, and all that remains is to choose suitable design sizes. For the quadrature rules in class 3 we naturally fix $m=N$, but for the regression designs in classes 1 and 2 we are free to choose any $m \geq N$. As described in \S\ref{PC}, we experiment with $m=N$ (\emph{uniquely-determined}) and $m=2N$ (\emph{twice over-determined}). This is the only difference between class 1 and class 2 designs. To give an example from Table \ref{designs}, for the \verb+adJULES+ experiment we have $n=4$, therefore a third-order polynomial chaos expansion ($p=3$) uses the following design sizes: $m = N = {4+3 \choose 3} = 35$ (class 1); $m = 2N = 2\times{4+3 \choose 3} = 70$ (class 2); $m = N = (3+1)^4 = 256$ (class 3). \\

For each design in Table \ref{designs} we build a polynomial chaos surrogate using the appropriate method to find the expansion coefficients. From this point on, all polynomial chaos surrogates fitted by the different non-intrusive approaches will be referred to by the acronym PC for simplicity. For all design classes and sizes in Table \ref{designs}, separate Gaussian process emulators with squared exponential and Mat\'{e}rn covariance functions are also built; in our experiments we refer to these using the acronyms SE GP and M GP respectively. The rationale of choosing our experimental designs provides a structure to our experiments, fixes the computational budget and ensures that we compare surrogates on a level playing field. The Sobol sequences in class 1 and 2 designs are generated randomly and kept the same for building each type of surrogate, so should not favour either surrogate method in any way. It would be beneficial to examine the effect of design choice through repeating the experiments with other randomly generated Sobol sequences, however this is not feasible because the simulators are expensive. Furthermore, sequentially or adaptively building the surrogates would no doubt lead to an improved performance in either case, but for a simple comparison we choose to fix the design a priori. \\

Since the sizes of the largest tensor grid designs in each experiment are much bigger than any of the Sobol sequence designs, we also use Sobol sequences of size $625$ and $1024$ in the \verb+adJULES+ and \verb+VEGACONTROL+ experiments respectively. This is to test whether there is any difference in the choice of experimental design when we can afford a large number of runs. Note also that a general rule of thumb in the Gaussian process literature is to use a design size ten times that of the input dimension ($m = 10n$) \citep{Loeppky:etal:2009}. Since $n=4$ and $5$ for the \verb+adJULES+ and \verb+VEGACONTROL+ simulators respectively, many of our design sizes are smaller than this suggestion. Hence, we are also testing the accuracy of surrogates built on small design sizes in our experiments. \\

To recap, in our experiments PC, SE GP and M GP surrogates are built on all the designs in Table \ref{designs}. For each of the \verb+adJULES+ and \verb+VEGACONTROL+ experiments we also have a independent $1000$-point Latin Hypercube design to be used for validation. Simulator output at the validation design points are compared with predictions from each of the surrogate methods. For each surrogate method and design, we compute the RMSE, mean, standard deviation, exceedance probabilities and PDF as described in \S\ref{Validation}. Surrogate mean, standard deviation, exceedance probability and PDF estimates are compared to corresponding simulator quantities evaluated from the validation design, whereas the RMSE is preferred to be a low as possible. As outlined in \S\ref{Introduction} and \S\ref{Validation} we wish to investigate how the accuracy of each surrogate method --- measured using the validation metrics --- is affected by changes in design size and class. In the process, we aim to identify scenarios in which polynomial chaos and Gaussian process methods might outperform one another.

\begin{table}
\centering
\footnotesize
\caption{Summary of the designs used in the adJULES and VEGACONTROL experiments. The design size $m$ increases to allow polynomial chaos expansions with larger truncation orders $p$. Note that a fourth class 3 design was not implemented for the VEGACONTROL experiment due to computational restrictions. Sobol sequence designs of size $625$ and $1024$ were also used for the adJULES and VEGACONTROL experiments respectively for comparison to the largest tensor grid designs. }
\begin{tabular}{|c|c|c|c|c|}
\hline
Class & Type & $p$ & $m_{\mathrm{adJULES}}$ & $m_{\mathrm{VEGACONTROL}}$ \\ \hline
\multirow{4}{*}{Class 1} & \multirow{4}{*}{Sobol sequence} & 1 & 5 & 6 \\
 & & 2 & 15 & 21 \\
 & & 3 & 35 & 56 \\
 & & 4 & 70 & 126 \\ \hline
\multirow{4}{*}{Class 2} & \multirow{4}{*}{Sobol sequence} & 1 & 10 & 12 \\
 & & 2 & 30 & 42 \\
 & & 3 & 70 & 112 \\
 & & 4 & 140 & 252 \\ \hline
\multirow{4}{*}{Class 3} & \multirow{4}{*}{Tensor grid} & 1 & 16 & 32 \\
 & & 2 & 81 & 243 \\
 & & 3 & 256 & 1024 \\
 & & 4 & 625 & --- \\ \hline
\end{tabular}
\label{designs}
\end{table}

\subsection{adJULES}
\label{adJULES}

\begin{figure}
\centering
\caption{Visualisation of the adJULES simulator output (cost function) as a function of pairwise combinations of each of the $4$ input parameters. The smoothing of the simulator output is performed using the surrogate with the lowest root mean square error metric in the adJULES experiments (see Figure \ref{JULESplot}), which in this case is the polynomial chaos expansion built on the largest class 3 design (see Table \ref{designs}).}
\label{JULESsim}
\epsfysize=115mm \epsfbox{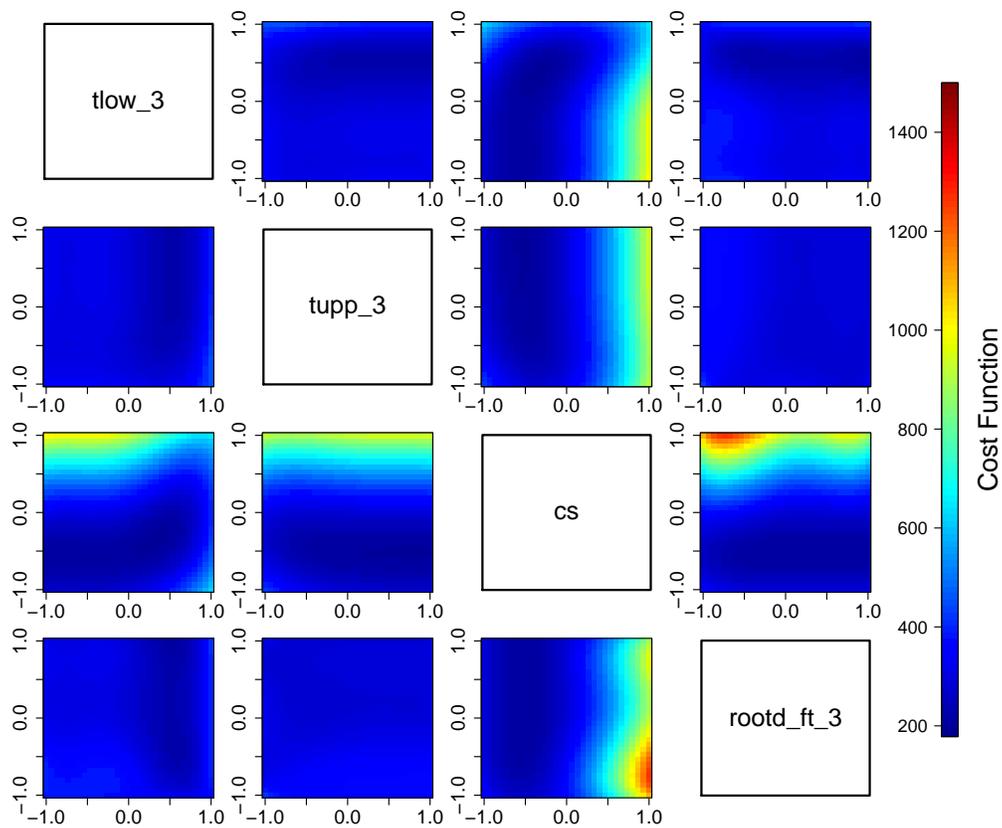}
\end{figure}

The Joint UK Land Environment Simulator (\verb+JULES+) \citep{JULES1,JULES2} is a simulator which models the interactions between the land and the atmosphere. It is the land surface component currently used in the UK Met Office Unified Model. \verb+JULES+ uses a number of  meteorological drivers and vegetation processes (e.g., photosynthesis, soil microbial activity) to model radiation, heat, water and carbon fluxes. Presently, observed time-series of these fluxes cannot be incorporated into the \verb+JULES+ framework. A new system, \verb+adJULES+, has been developed to provide this functionality and also comprises an adjoint model for parameter estimation or optimisation studies. The current implementation of \verb+adJULES+ contains $9$ input parameters detailing various plant properties and a single output, a cost function for optimisation purposes. The gradient of the output with respect to each of the input parameters is also available, but for simplicity we do not make use of this in building the surrogates. Expert elicitation prior to the experiment led us to focus on $4$ of the most important input parameters; these were \verb+t_low+, \verb+t_upp+, \verb+cs+ and \verb+rootd_ft+, which represent lower and upper temperatures for photosynthesis ($\degree \mathrm{C}$), carbon content of soil ($\mathrm{kg} \, \mathrm{C} \, \mathrm{m}^{-2}$) and root depth (m) respectively. To build the surrogates, we first transform the parameters to be distributed on $[-1,1]$, so that $\mathcal{D} \subseteq [-1,1]^4$, and assume a uniform distribution across this space. We keep the remaining $5$ parameters at their nominal values when executing the \verb+adJULES+ simulator. We consider the \verb+adJULES+ simulator to be a black-box and aim to best represent the simulator output (cost function) across the space of the $4$ input parameters described above. We do this by building (non-intrusively) PC, SE GP and M GP surrogates on the designs in Table \ref{designs} and seek to compare their accuracy, as measured by the validation metrics presented in Section \ref{Validation}, with changes in design type and class. A visualisation of the \verb+adJULES+ simulator output (cost function) as a function of pairwise combinations of each of the $4$ input parameters can be seen in Figure \ref{JULESsim}. \\

\begin{figure}
\centering
\caption{Root mean square error, mean and standard deviation validation results for the adJULES simulator. Design classes are in columns and validation metrics are in rows. Polynomial chaos, squared exponential and Mat\'{e}rn Gaussian process surrogate models are shown as red triangles, purple crosses and blue circles respectively. The points have been been jittered slightly for clarity. Gaussian process emulators also have $95\%$ confidence intervals about their estimates (solid lines). Simulator mean and standard deviation (solid black lines) are also shown with $95\%$ confidence intervals (dashed black lines).}
\label{JULESplot}
\epsfysize=105mm \epsfbox{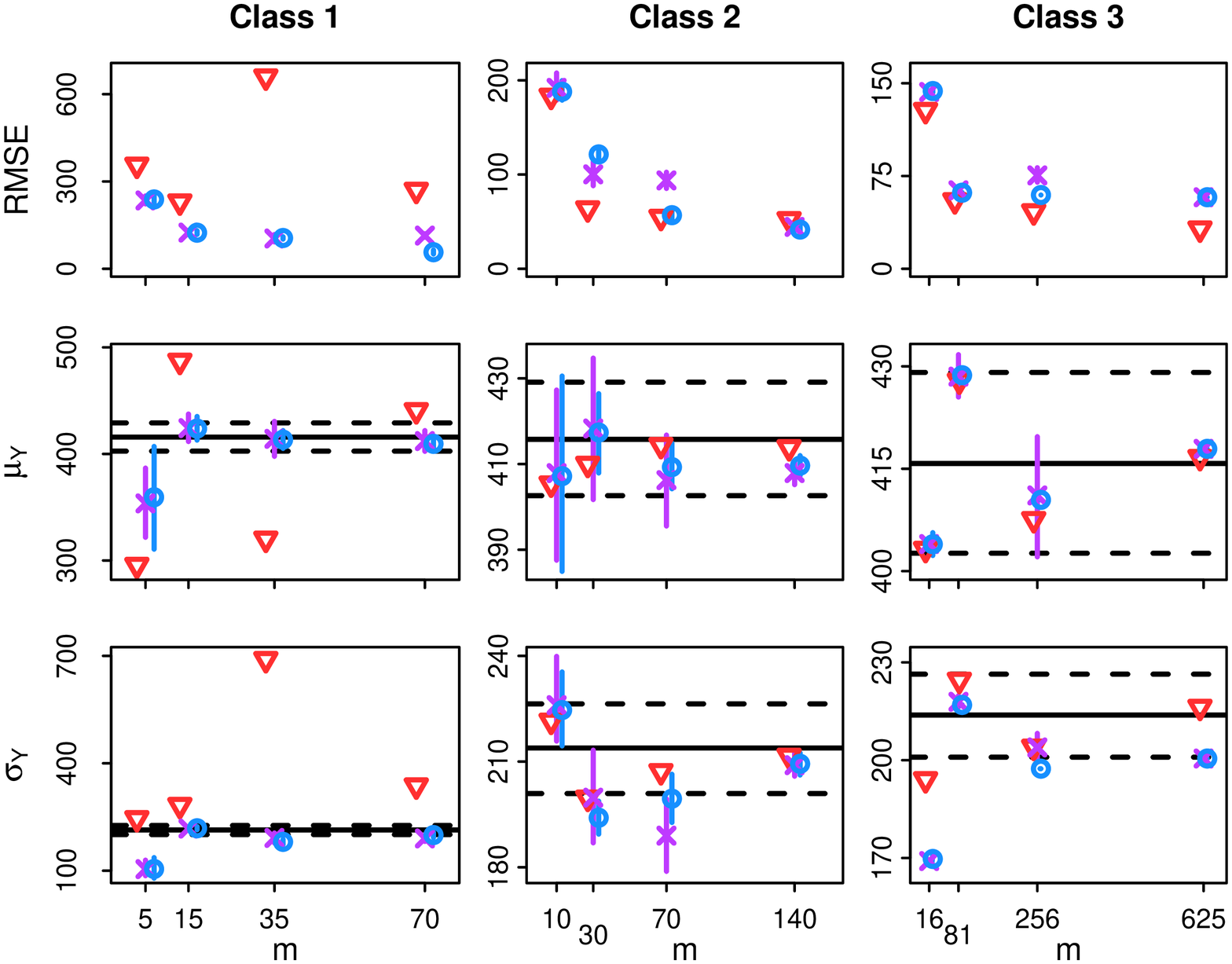}
\end{figure}

\begin{figure}
\centering
\caption{Exceedance probability validation results for the adJULES simulator. Design classes are in columns and validation metrics are in rows. Polynomial chaos, squared exponential and Mat\'{e}rn Gaussian process surrogate models are shown as red triangles, purple crosses and blue circles respectively. The points have been been jittered slightly for clarity. Gaussian process emulators also have $95\%$ confidence intervals about their estimates (solid lines). Simulator exceedance probabilities (solid black lines) are also shown with $95\%$ confidence intervals (dashed black lines).}
\label{JULESexceed}
\epsfysize=90mm \epsfbox{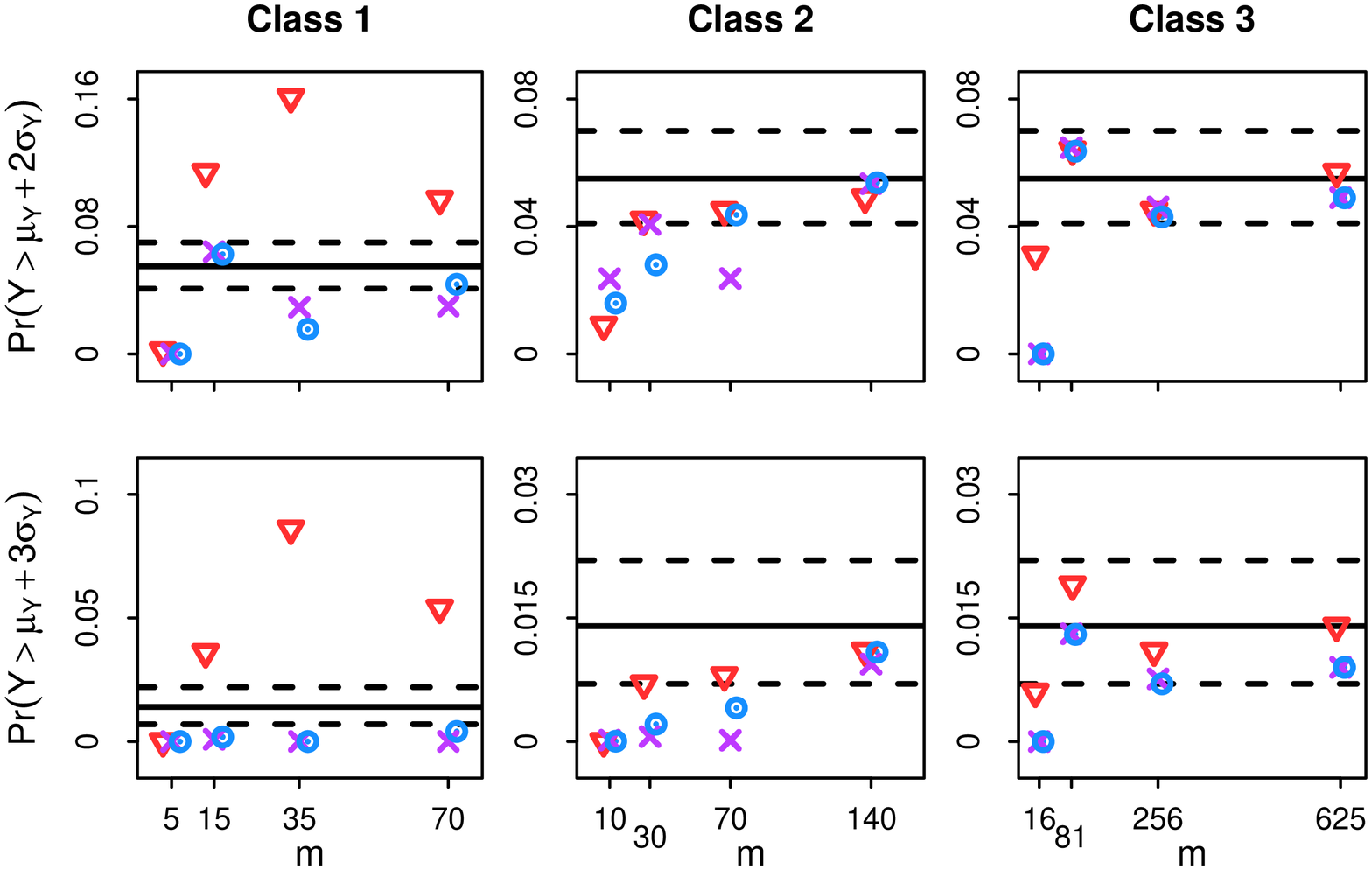}
\end{figure}

\begin{figure}
\centering
\caption{Probability density function (PDF) estimation validation results for the adJULES simulator. Design classes are in columns and design size used to build the surrogates increases further down the rows (see Table \ref{designs} for exact design size in each panel). Polynomial chaos, squared exponential and Mat\'{e}rn Gaussian process surrogate model PDF estimates are shown as red, purple and blue lines respectively, while the simulator output PDF is shown as a black line. Gaussian process emulators also have $95\%$ confidence intervals about their estimates.}
\label{JULESdensity}
\epsfysize=110mm \epsfbox{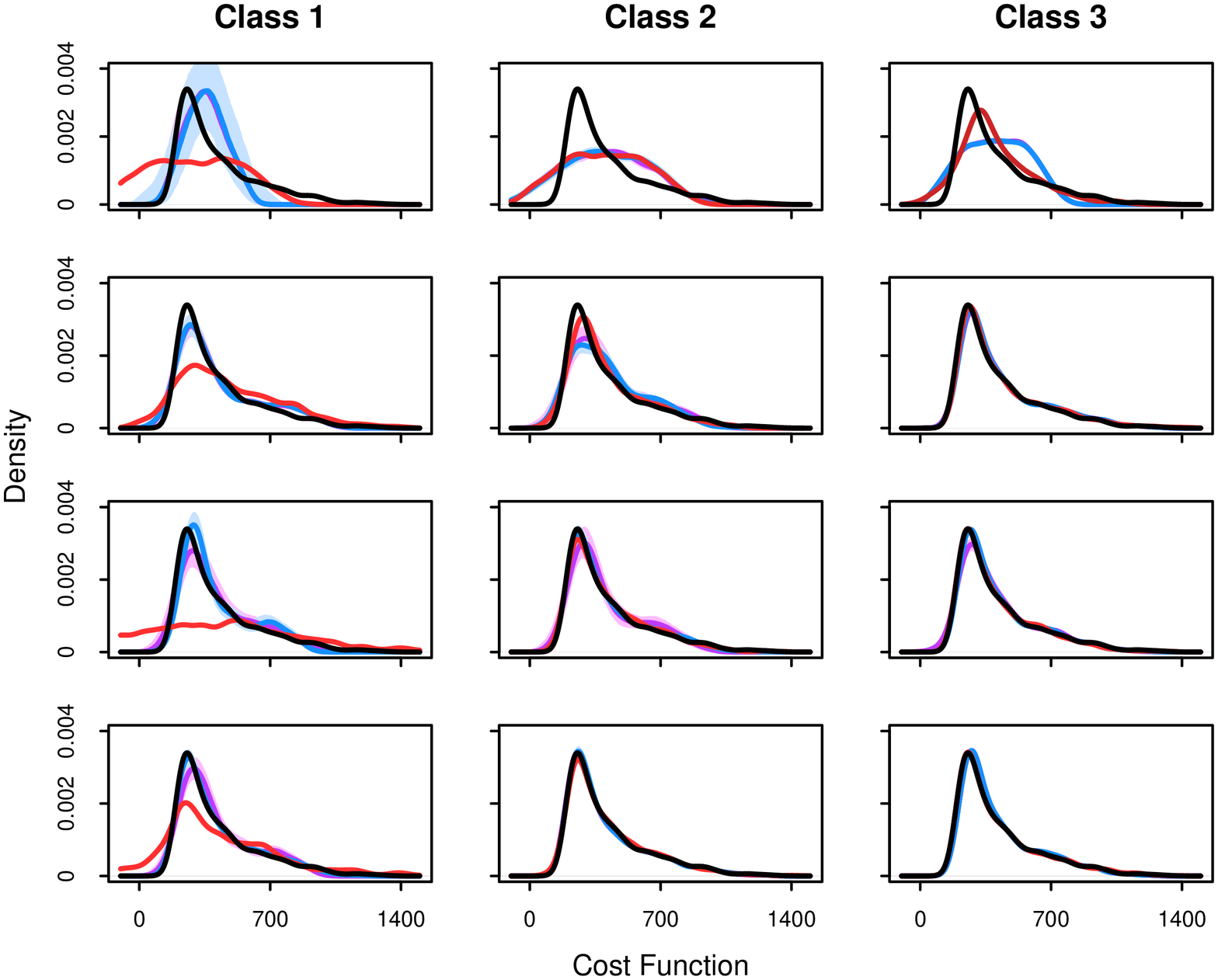}
\end{figure}

Validation results from fitting PC, SE GP and M GP surrogates to the \verb+adJULES+ simulator are presented in Figures \ref{JULESplot}, \ref{JULESexceed} and \ref{JULESdensity}. Firstly examining Figure \ref{JULESplot}, we note that for class 1 designs both GP emulator types outperform PC for the RMSE, mean and standard deviation metrics. The accuracy of PC surrogates does not necessarily improve with design size and consequently the addition of more expansion terms, suggesting that uniquely-determined regression approach is unstable. Particularly we draw attention to the outliers for PC in the third design of this class (size $35$). For this design class there is little difference between the two GP emulators, and they are the preferred method. For class 2 designs, twice over-determined regression is used to determine the expansion coefficients, leading to more interesting results. When analysing RMSE there is little to choose between GP and PC approaches. We note that PC has a faster initial reduction in error but does not continue to improve. Both GP methods have a slower reduction in error but may begin to outperform PC methods for larger designs in this class. Concerning the simulator mean, all surrogate estimates are within our calculated confidence intervals. The uncertainty on the GP estimates can be seen to reduce with larger design size. For the simulator standard deviation, PC and GP estimates are mostly similar although we note that PC estimates tend to be closer to the `true' value. However, uncertainty on the GP estimates themselves does mean it is difficult to choose a favoured surrogate approach here. For class 3 designs, we observe that PC methods are consistently more accurate than GP emulators for the RMSE metric so they are preferred here. Similarly to results from class 2 designs, when estimating the simulator mean the surrogate approaches are very consistent with one another and are all within the confidence intervals. Some preference for PC methods can be found in the case of standard deviation, as GP estimates are unstable to increases in design size whereas PC estimates gradually become more accurate. \\

Secondly, we analyse results from the exceedance probability metrics shown in Figure \ref{JULESexceed}. For class 1 designs it is clear that GP methods outperform PC surrogates for all design sizes and both exceedance probabilities. However we note that while GP estimates are comparatively better and preferred here, they are still not particularly accurate. Furthermore, the uncertainty attached to GP estimates are always small meaning they are confident about their inaccurate predictions. Clearly for these small design sizes it is different to get an accurate estimate of such small probabilities --- no surrogate method can get within the confidence intervals calculated for the probability of exceeding three standard deviations above the mean. For class 2 designs we see a great improvement for PC methods, and they arguably give the most accurate estimates for both exceedance probabilities and most design sizes in the class, falling within the calculated simulator confidence intervals more often. One can also observe more stable results here compared to class 1, with accuracy generally improving with design size for all surrogate methods. For class 3 designs, results are similar between PC and GP methods and their estimates fall within the confidence intervals, with the exception of the smallest design size considered. It could perhaps be said here that PC methods do better at estimating the probability of exceeding three standard deviations above the simulator mean. \\

Finally, we examine the capability of the surrogates at estimating the simulator output PDF, shown in Figure \ref{JULESdensity}. As observed with the previous validation metrics, a clear preference for GP methods can be found for class 1 designs, where PC methods struggle to estimate the PDF accurately in all cases. For class 2 and 3 designs, it is difficult to observe many differences between the methods, although PC can be seen to have a slight advantage for individual designs (for example the second and third designs in class 2, and the first design in class 3). For large design sizes in these classes all methods can estimate the PDF to a high degree of accuracy. \\

In summary, GP methods are preferred for class 1 designs regardless of validation metric; RMSE results are similar for surrogates built on class 2 designs but PC is favoured for class 3 designs; all surrogate methods accurately estimate the simulator mean for class 2 and class 3 designs; PC is narrowly favoured for estimating the simulator standard deviation in all cases except for class 1 designs. For both the exceedance probabilities considered, neither method is accurate for class 1 designs but GP is preferred; PC is more accurate for class 2 designs; methods perform similarly for class 3 designs with a narrow preference for PC. Finally, density estimation results show that both methods give accurate results for large design sizes in classes 2 and 3. Some preferences can be found for PC in the smaller designs in these classes. \\

We also built PC, SE GP and M GP surrogates for the \verb+adJULES+ simulator on a Sobol sequence design of size $625$, to see if any improvement could be made on the validation results from the largest tensor grid design. RMSE, mean and standard deviation validation metrics for surrogates built on tensor grid and Sobol sequence designs of size $625$ are shown in the top panels of Figure \ref{sobol}. We observe that with the use of a large Sobol sequence design, both GP emulators become more accurate in terms of RMSE and now have similar results to PC. While the simulator mean estimates still remain accurate for all surrogate methods with the design change, for simulator standard deviation the GP emulator estimates are now within the confidence intervals. Validation metrics for PC do not worsen with the change to Sobol sequence designs and remain highly accurate. We expect these results to carry over to the exceedance probability and PDF estimation validation metrics, but do not consider these here for simplicity.

\subsection{VEGACONTROL}
\label{VEGACONTROL}

\verb+VEGACONTROL+ is a Matlab C-coded simulator used in industry to prepare and validate the VEGA flight management and control system in the atmospheric and exo-atmospheric flight phases \citep{Mujumdar:etal:2015}. It is a non-linear, six degrees of freedom industrial simulation model of the VEGA launch vehicle, which is a recent European multi-payload launch vehicle developed by European Launch Vehicle under European Space Agency (ESA) responsibility{\footnote{The simulation environment has been used as a benchmark model for the ESA research activity ``Robust Flight Control System Design Verification and Validation Framework'' (ESA AO/1-6322/09/NL/JK).}}. The flight phase with altitude between 30m and 60km is considered in the present study. The simulator equations for motion include force and drag components depending on Mach and angle of attack, kinematic coupling in all axes, and a non-linear model of the electro-mechanical actuator dynamics with associated backlash and delays. The mathematical model for QUASAR Inertial Sensor Unit with its noise and bias characteristics and the propulsion model of the P80 solid propulsion system with validated thrust curves that include thrust oscillation effects to assess proper execution of the separation dynamics further constitute to the complexity of the simulator. \\

A high fidelity structural flexible mode model describing the launcher deformation is included to assess proper filtering and stability properties.  The atmosphere model includes also a set of measured sizing wind-gust input models representative for the Kourou launch site. The launcher dynamics are driven by the FPSA ADA/C-flight code reflecting the flight management system for the time line sequence command and execution of associated guidance navigation and control system for thrust vector control (TVC) and roll and attitude control (RACS) and other support functions such as acceleration threshold detection and pyro-valve command for stage separation.\\

\begin{figure}
\centering
\caption{Visualisation of the VEGACONTROL simulator output $\mathrm{max}(Q_\alpha)$ as a function of pairwise combinations of each of the $5$ input parameters. The smoothing of the simulator output is performed using the surrogate with the lowest root mean square error metric in the VEGACONTROL experiments (see Figure \ref{VEGAplot}), which in this case is the Gaussian process emulator with Mat\'{e}rn covariance function built on the largest class 2 design (see Table \ref{designs}).} 
\label{VEGAsim}
\epsfysize=115mm \epsfbox{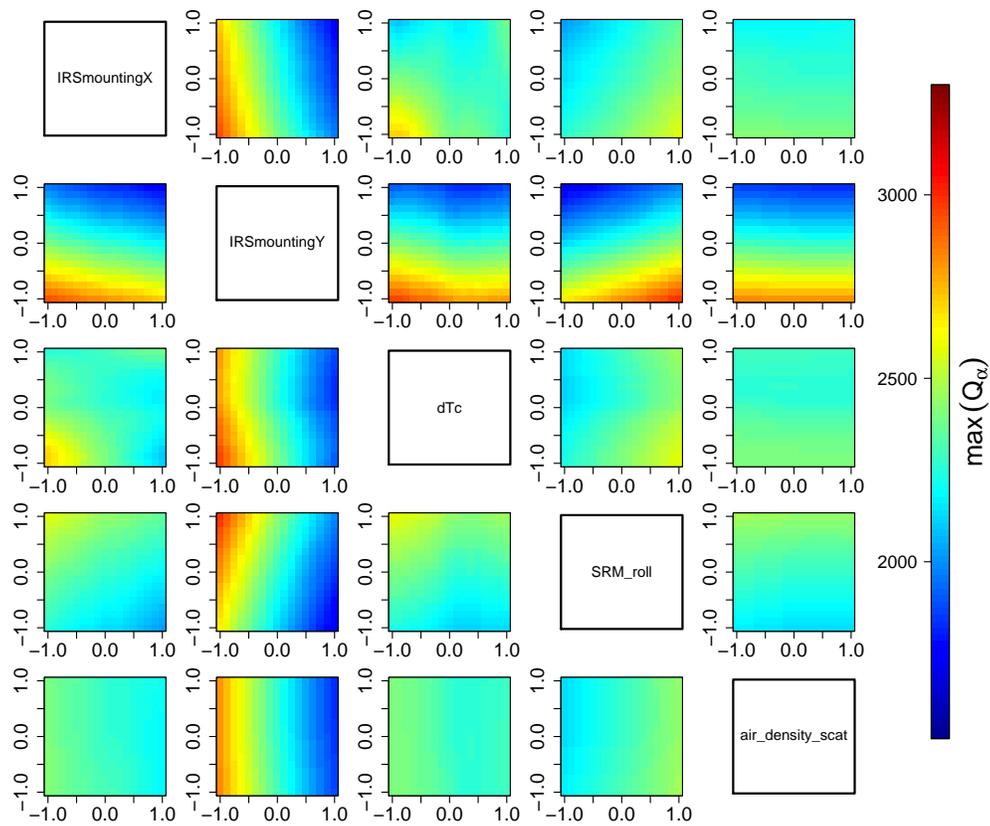}
\end{figure}

The overall validation and verification criteria considered in this research activity represent the TVC and RACS technical requirements for the atmospheric phase of flight. Load requirement deviations ($Q_\alpha$) must be limited via the product of dynamic pressure $Q$ and angle of attack $\alpha$ over the entire P80 flight Mach range. The problem addressed in this paper attempts to assess the requirement on $\mathrm{max}(Q_\alpha)$ over a reduced parameter combination subset that dominates the requirement degradation. The current implementation of \verb+VEGACONTROL+ has $83$ input parameters and $37$ outputs. Expert elicitation prior to the experiment led us to focus on $5$ of the most influential input parameters; these were \verb+IRSmountingX+, \verb+IRSmountingY+, \verb+dTc+, \verb+SRM_roll+ and \verb+air_density_scat+. These parameters correspond to the Indian Remote Sensing (IRS) mounting error with respect to the X and Y body axes respectively, scattering on time burn, scattering on P80 roll degree and atmospheric density respectively. As already mentioned we focus on a single output, namely $\mathrm{max}(Q_\alpha)$, which is the maximum of the aerodynamic load. The gradient of this output with respect to the input parameters is not available. As with the \verb+adJULES+ experiment we transform the input parameters to be uniformly distributed on $[-1,1]$, so that here $\mathcal{D} \subseteq [-1,1]^5$, and keep the remaining input parameters at their nominal values when executing the simulator. Once again, we consider the \verb+VEGACONTROL+ simulator to be a black-box and aim to best represent the simulator output ($\mathrm{max}(Q_\alpha)$) across the space of the $5$ input parameters described above. We do this by building (non-intrusively) PC, SE GP and M GP surrogates on the designs in Table \ref{designs} and seek to compare their accuracy, as measured by the validation metrics presented in Section \ref{Validation}, with changes in design type and class. A visualisation of the \verb+VEGACONTROL+ simulator output $\mathrm{max}(Q_\alpha)$ as a function of pairwise combinations of each of the $5$ input parameters can be seen in Figure \ref{VEGAsim}. \\

\begin{figure}
\centering
\caption{Root mean square error, mean and standard deviation validation results for the VEGACONTROL simulator. Design classes are in columns and validation metrics are in rows. Polynomial chaos, squared exponential and Mat\'{e}rn Gaussian process surrogate models are shown as red triangles, purple crosses and blue circles respectively. The points have been been jittered slightly for clarity. Gaussian process emulators also have $95\%$ confidence intervals about their estimates (solid lines). Simulator mean and standard deviation (solid black lines) are also shown with $95\%$ confidence intervals (dashed black lines).}
\label{VEGAplot}
\epsfysize=105mm \epsfbox{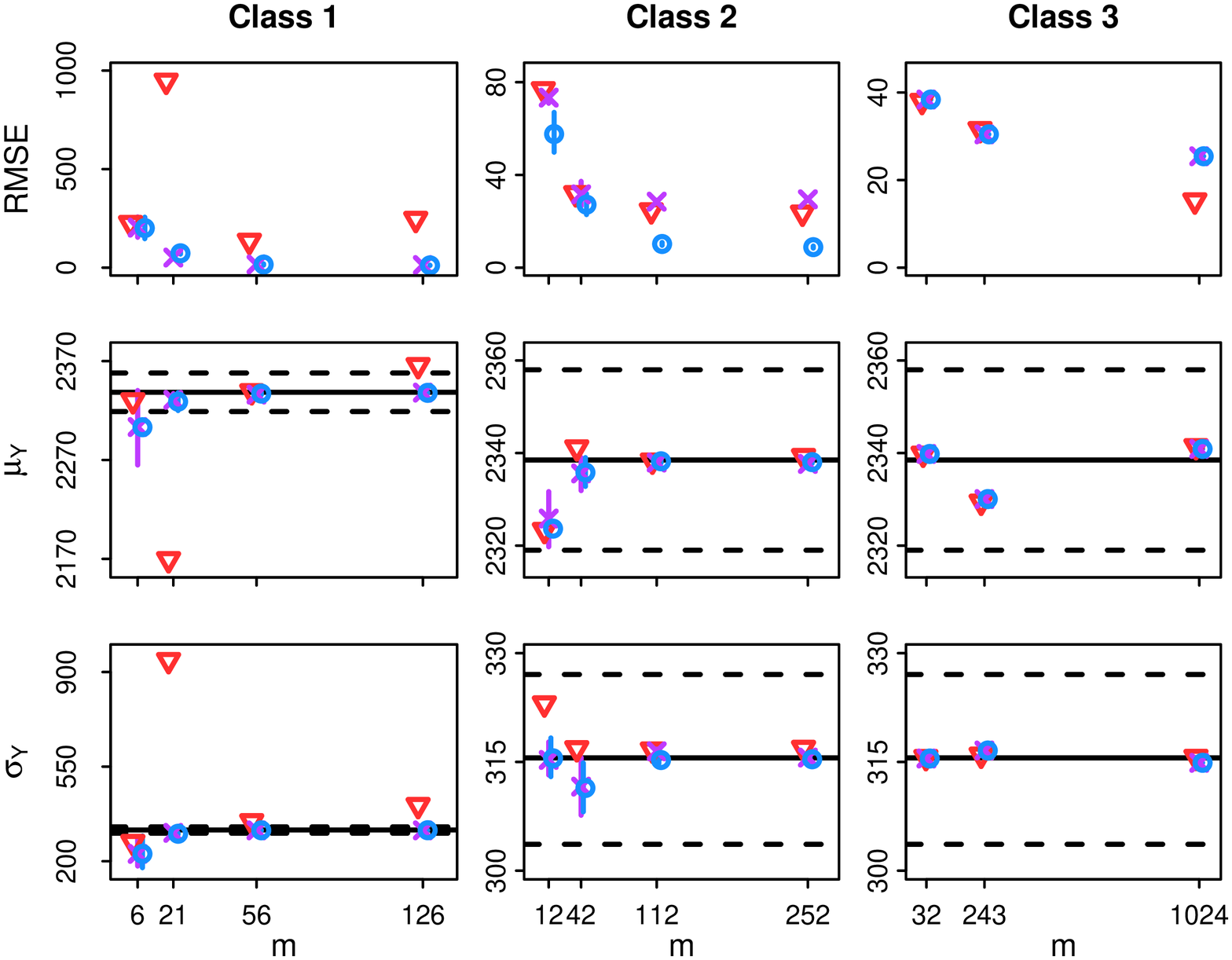}
\end{figure}

\begin{figure}
\centering
\caption{Exceedance probability validation results for the VEGACONTROL simulator. Design classes are in columns and validation metrics are in rows. Polynomial chaos, squared exponential and Mat\'{e}rn Gaussian process surrogate models are shown as red triangles, purple crosses and blue circles respectively. The points have been been jittered slightly for clarity. Gaussian process emulators also have $95\%$ confidence intervals about their estimates (solid lines). Simulator exceedance probabilities (solid black lines) are also shown with $95\%$ confidence intervals (dashed black lines).}
\label{VEGAexceed}
\epsfysize=90mm \epsfbox{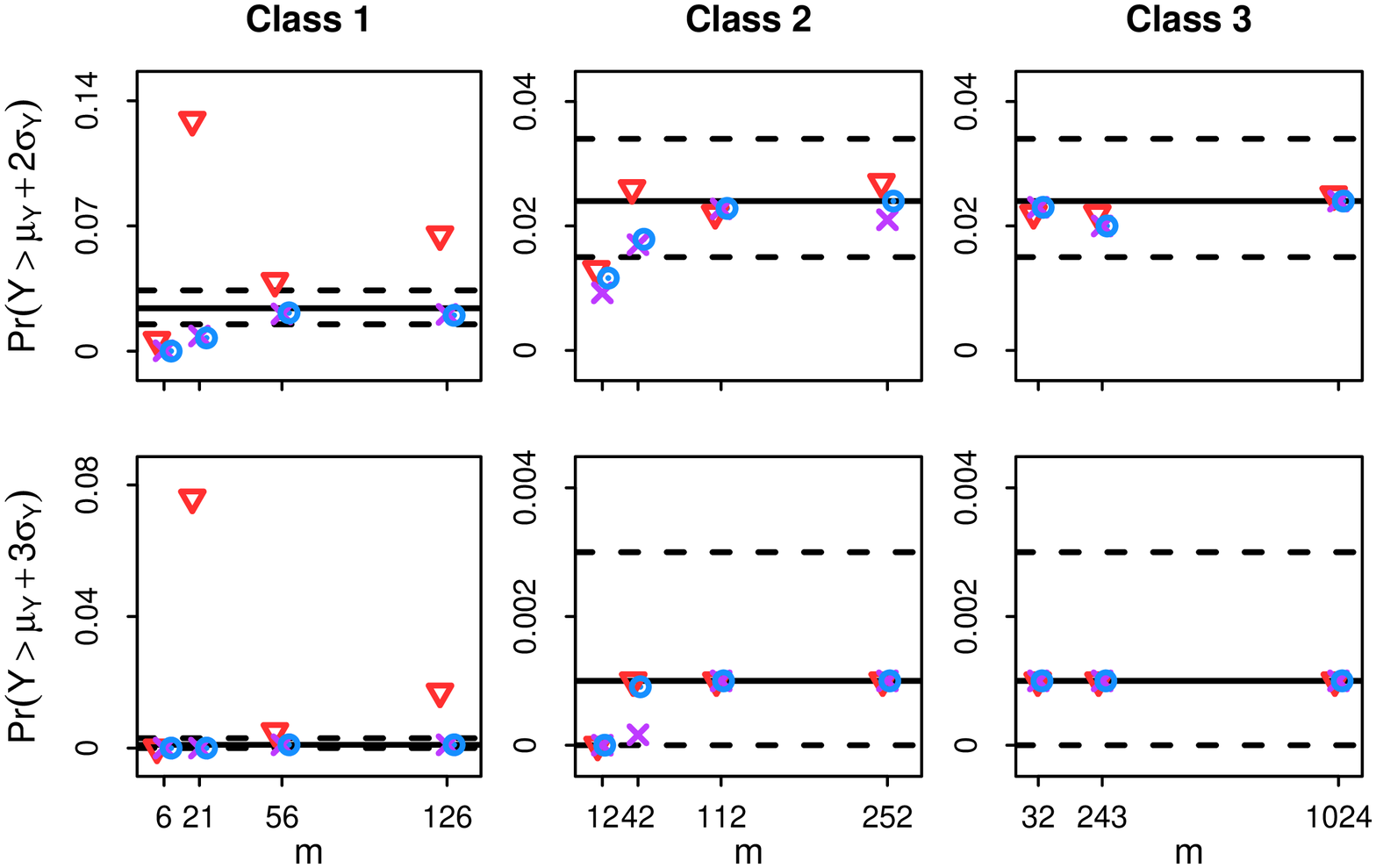}
\end{figure}

\begin{figure}
\centering
\caption{Probability density function (PDF) estimation validation results for the VEGACONTROL simulator. Design classes are in columns and design size used to build the surrogates increases further down the rows (see Table \ref{designs} for exact design size in each panel). Polynomial chaos, squared exponential and Mat\'{e}rn Gaussian process surrogate model PDF estimates are shown as red, purple and blue lines respectively, while the simulator output PDF is shown as a black line. Gaussian process emulators also have $95\%$ confidence intervals about their estimates.}
\label{VEGAdensity}
\epsfysize=110mm \epsfbox{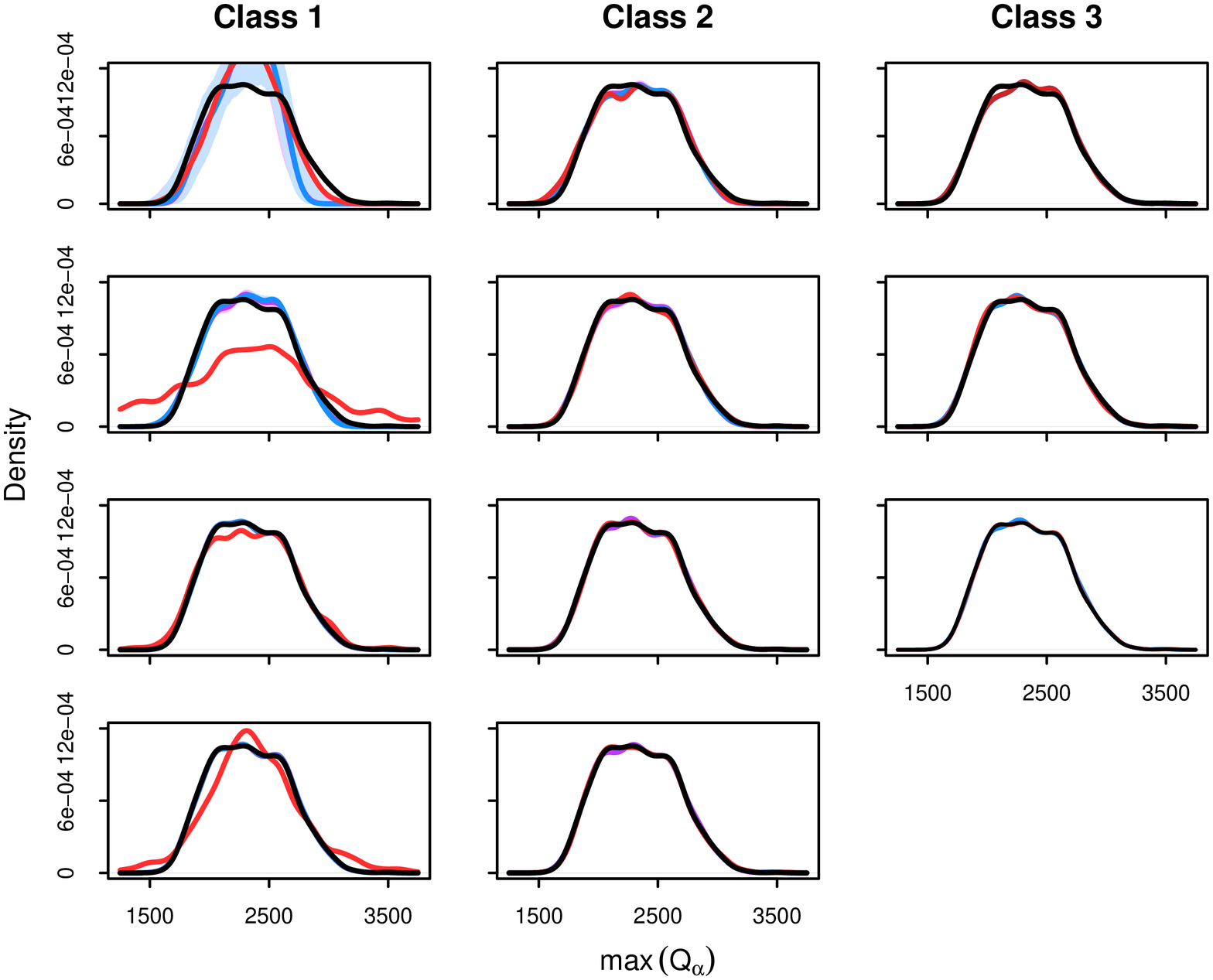}
\end{figure}

Validation results for the \verb+VEGACONTROL+ simulator are shown in Figures \ref{VEGAplot}, \ref{VEGAexceed} and \ref{VEGAdensity}. Looking at Figure \ref{VEGAplot}, once again we observe that for class 1 designs, both the GP emulators outperform PC for the RMSE, mean and standard deviation metrics. In particular, the PC surrogate built on the second design of this class (size $21$) has a large error for all validation metrics. For class 2 designs, all the surrogate approaches exhibit similar accuracy. This is particularly the case when estimating simulator mean and standard deviation; we find that surrogate estimates are similar to one another and are also all within the calculated confidence intervals. In this case, none of the surrogate methods would be preferred. However, when looking at RMSE results there is some evidence for favouring GP emulation over PC. While PC and SE GP surrogates have similar accuracy for all design sizes, the M GP surrogate consistently does better. This is one of the few cases where there is substantial difference in choice of covariance function, and perhaps the Mat\'{e}rn covariance function is particularly suited to output from the \verb+VEGACONTROL+ simulator. The converse is found for class 3 designs. While again the surrogate estimates of simulator mean and standard deviation are almost identical and of high quality, there is some evidence to suggest that PC methods perform better in terms of RMSE. However, this can only be observed for the largest design in this class; for the two smaller designs all the surrogate approaches have a similar RMSE. \\

Examining the exceedance probabilities in Figure \ref{VEGAexceed}, we see that all surrogate methods perform better than in the \verb+adJULES+ experiment. We attribute this to the fact that the simulator output PDF is more symmetrical and has a shorter tail, although perhaps the \verb+VEGACONTROL+ simulator is better behaved in the tails. Once again for class 1 designs the PC estimates are inaccurate and unstable with increases in design size, and GP methods are favoured. For class 2 and 3 designs, estimates from GP and PC surrogates are similar and of high quality for both exceedance probabilities and all design sizes considered. In this case, none of the surrogate methods would be preferred. \\

Lastly we consider PDF estimation, shown in Figure \ref{VEGAdensity}. A clear preference for GP methods can again be found for class 1 designs, where PC methods struggle to estimate the PDF accurately in the majority of cases. For class 2 and 3 designs, all surrogate PDF estimates are of high quality and there is little to choose between the methods. \\ 

In summary, GP methods are preferred for class 1 designs regardless of validation metric; mean, standard deviation, exceedance probability and PDF estimates are high quality and similar between surrogate approaches for class 2 and 3 designs; for lower RMSE the M GP emulator is preferred for class 2 designs, whereas PC is favoured for class 3 designs. In general, the results from the \verb+adJULES+ experiment nicely carry over to the \verb+VEGACONTROL+ simulator apart from the fact that all surrogates tended to calculate the validation metrics to a higher degree of accuracy here. \\

\begin{figure}
\centering
\caption{Difference in root mean square error, mean and standard deviation estimates when using large Sobol sequences instead of tensor grids. Sobol sequences of size $625$ and $1024$ were used in the adJULES (top panels) and VEGACONTROL (bottom panels) experiments respectively to correspond with the largest tensor grid designs. Validation diagnostics from the tensor grid designs are shown in black. Diagnostics from the Sobol sequence designs are shown in colour: polynomial chaos, squared exponential and Mat\'{e}rn Gaussian process surrogate models are shown as red triangles, purple cross and blue circles respectively.}
\label{sobol}
\epsfysize=90mm \epsfbox{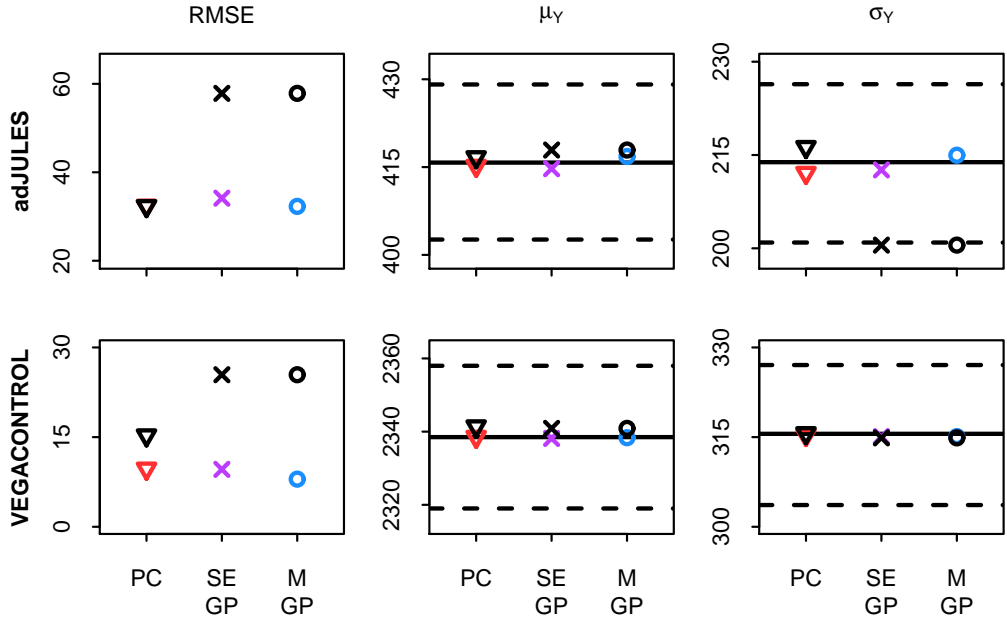}
\end{figure}

We also built PC, SE GP and M GP surrogates for the \verb+VEGACONTROL+ simulator on a Sobol sequence design of size $1024$, to see if any improvement could be made on the validation results from the largest tensor grid design. RMSE, mean and standard deviation validation metrics for surrogates built on tensor grid and Sobol sequence designs of size $1024$ are shown in the bottom panels of Figure \ref{sobol}. We note that with the use of a large Sobol sequence design, both GP emulators become more accurate in terms of RMSE. Some improvement can also be observed for PC, and now all surrogate methods show similar accuracy. Estimates of the simulator mean and standard deviation remain accurate and are robust under the change of design. Again we expect these results to carry over to the exceedance probability and PDF estimation validation metrics, but do not consider these here for simplicity.       

\section{Discussion}
\label{Discussion}

In this paper we have compared two popular surrogate methods for uncertainty quantification in computer experiments, polynomial chaos and Gaussian process emulation. Based on results from fitting surrogates to the output of two simulators, we have found that one approach does not unanimously outperform the other but that the best method depends on the modelling goals of the practitioner and on the type of experimental design. If one requires an accurate surrogate across the design space (which we measure using RMSE), we generally observed that Gaussian process emulators are preferred for Sobol sequence type designs, whereas polynomial chaos was more accurate for tensor grid designs. We have two possible reasons for why this should be. Firstly, research has found (for example, \citet{Urban:Fricker:2010}) that Gaussian process emulators are more accurate when built on less structured designs than tensor grids. This is mainly because if we project a multidimensional tensor grid design onto a single dimension, some design points are repeated and do not provide extra information. This does not happen for Sobol sequences or Latin Hypercubes, and the extra data in each dimension results in better estimation of Gaussian process parameters. Secondly, the regression approach used for polynomial chaos on the Sobol sequence designs leads to surrogates which do not necessarily interpolate the simulator output, but the quadrature approach on tensor grid designs does. This has some parallels with the use of a nugget parameter in Gaussian processes, where the interpolation property of the emulator is relaxed. If instead one simply wants to extract some properties of the simulator output (for example, the mean or standard deviation) there is a slightly different focus. When estimating the simulator mean, results showed that both Gaussian process emulation and polynomial chaos gave very similar and highly accurate results even for small design sizes. In terms of the standard deviation, a similar story was seen for the \verb+VEGACONTROL+ simulator but a narrow preference for polynomial chaos could be found in the \verb+adJULES+ experiment. When estimating the probability of exceeding two or three standard deviations above the simulator mean, while results for each surrogate were similar for the \verb+VEGACONTROL+ experiment, some small preferences for PC could be seen for the \verb+adJULES+ simulator (particularly for class 2 designs). In terms of estimating the simulator PDF, there was little to choose between the surrogate approaches. We also found that if one is restricted to small design sizes the uniquely-determined regression technique for estimating the polynomial chaos coefficients should be avoided, agreeing with \citet{Hosder:etal:2007} who suggest twice over-determined regression gives more stable results. In this case the practitioner should ensure enough design points for at least twice over-determined regression, or use Gaussian process emulation instead. Alternatively, sparse polynomial chaos techniques, which use regularisation to reduce the number of terms in the expansion (for instance, using LASSO or least-angle regression), may also lead to improved stability in this case. Lastly, given that there was some evidence of a performance gap between polynomial chaos and Gaussian process emulators for the largest tensor grid designs, we investigated how the validation metrics changed with a Sobol sequence design of the same size. We found that the performance gap could be closed completely if one used a Sobol sequence design instead. In particular, the RMSE for the Gaussian process emulators massively reduced to be consistent with polynomial chaos approaches. We reiterate that the choice of design is very important when building a surrogate; if one can afford a reasonable amount of design points then a tensor grid design would be suitable for polynomial chaos, whereas a Sobol sequence or Latin Hypercube design would be more appropriate for a Gaussian process emulator.  \\

Aside from the results of our experiments, we should also reflect on some generic advantages and disadvantages that polynomial chaos and Gaussian process emulation may have in relation to one another. Notably we draw upon some points concerned with the computational cost of building the surrogate, as well as their flexibility and practicality in different scenarios. In terms of the computational cost of building the surrogate, for polynomial chaos most of the work is done prior to the experiment in choosing the polynomial basis used in the expansion. Recall that this solely depends on the probability distributions chosen for the simulator input parameters. The subsequent non-intrusive computation of the expansion coefficients is generally cheap, requiring the evaluation of plug-in formulas (quadrature) or some simple matrix algebra (regression). Therefore if one simply wants to find the simulator mean or standard deviation, polynomial chaos methods would be the faster approach. Conversely, most of the computational cost comes from fitting the Gaussian process emulator itself, and this can be expensive (particularly for large design sizes). Problems can also arise when estimating the correlation lengths and the optimisation process for the parameters may have to be restarted. However, these problems usually have a simple fix (for example with the addition of a nugget parameter, see \citet{Andrianakis:Challenor:2012}). \\

Regarding the comparative flexibility of the surrogate approaches to adapt to different scenarios, it is clear that Gaussian process emulation has the edge. In all of our experiments we were restricted in our choice of design by the polynomial chaos method --- whether it be choosing Gaussian abscissae for the quadrature method or restricting the size of design to ensure a certain truncation order for the polynomial expansion. In this sense we gave a subtle advantage to polynomial chaos methods in our experiments. In contrast, there are really no restrictions on the type and size of design used to build a Gaussian process emulator; although naturally we adhere to principles such as a space-filling design, and require an adaquate design size for an accurate emulator. The fact that the Gaussian process emulator performed as well as polynomial chaos based on our design choices is testament to flexibility of the approach. We must also remark that although many simulators can be quickly and accurately modelled by a polynomial function, some are bound to exhibit more complicated behaviour. There is no doubt that Gaussian process emulators can deal with a wider range of simulator behaviours, and this can be done relatively easily with changes in the mean and covariance functions. To add more weight to this argument, consider the following two-dimensional function as a simulator,

\begin{equation}
y = \eta(x_1,x_2) = \exp(-x_1) \tanh(5 x_2) \qquad x_1,x_2 \in [-1,1] .
\label{nonlinear}
\end{equation}

\begin{figure}
\centering
\caption{Plot of the two-dimensional function in \eqref{nonlinear} for $x_1,x_2 \in [-1,1]$.}
\label{exampleplot}
\epsfysize=105mm \epsfbox{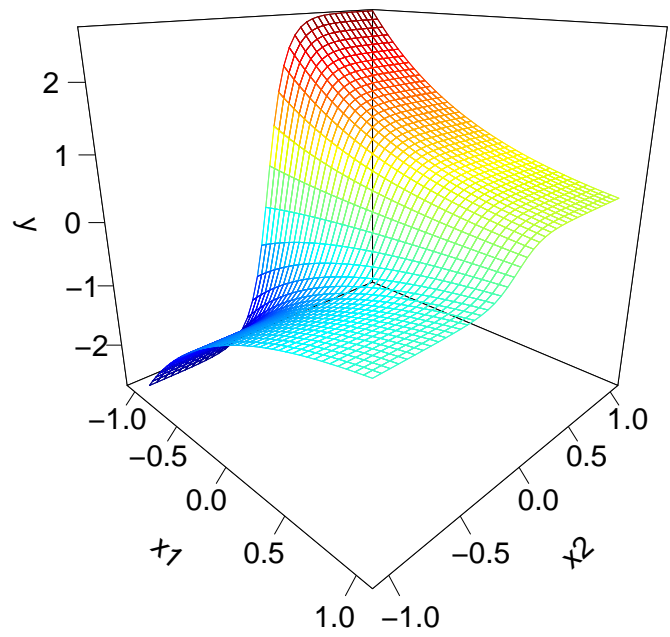}
\end{figure}

\begin{figure}
\centering
\caption{Root mean square error, mean and standard deviation validation results for the two-dimensional function in \eqref{nonlinear}. Design classes are in columns and validation metrics are in rows. Polynomial chaos, square exponential and Mat\'{e}rn Gaussian process surrogate models are shown as red triangles, purple crosses and blue circles respectively. The points have been been jittered slightly for clarity. Gaussian process emulators also have $95\%$ confidence intervals about their estimates (solid lines). Simulator mean and standard deviation (solid black lines) are also shown with $95\%$ confidence intervals (dashed black lines).}
\label{example}
\epsfysize=105mm \epsfbox{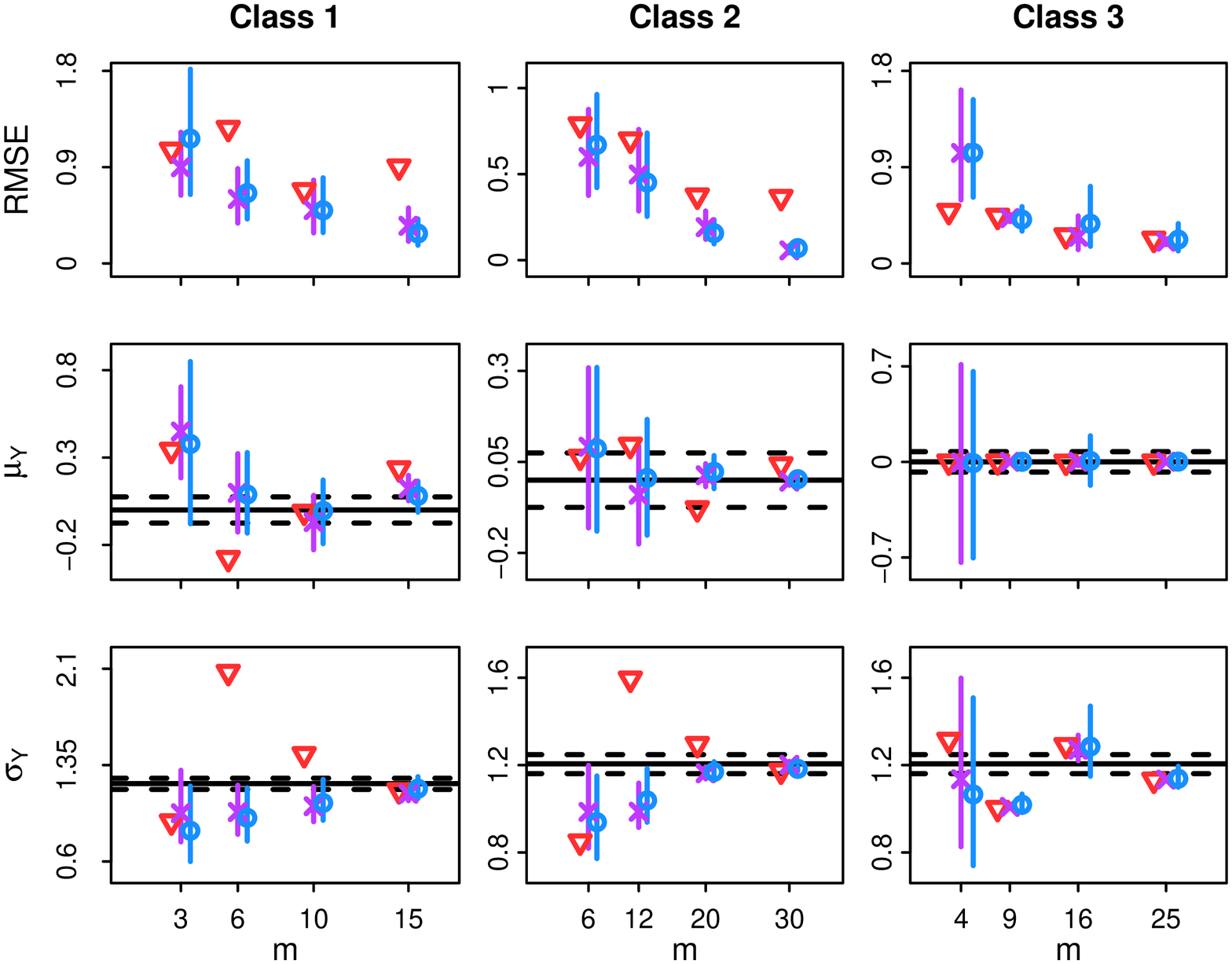}
\end{figure}

\begin{figure}
\centering
\caption{Exceedance probability validation results for the two-dimensional function in \eqref{nonlinear}. Design classes are in columns and validation metrics are in rows. Polynomial chaos, squared exponential and Mat\'{e}rn Gaussian process surrogate models are shown as red triangles, purple crosses and blue circles respectively. The points have been been jittered slightly for clarity. Gaussian process emulators also have $95\%$ confidence intervals about their estimates (solid lines). Simulator exceedance probabilities (solid black lines) are also shown with $95\%$ confidence intervals (dashed black lines).}
\label{exampleexceed}
\epsfysize=90mm \epsfbox{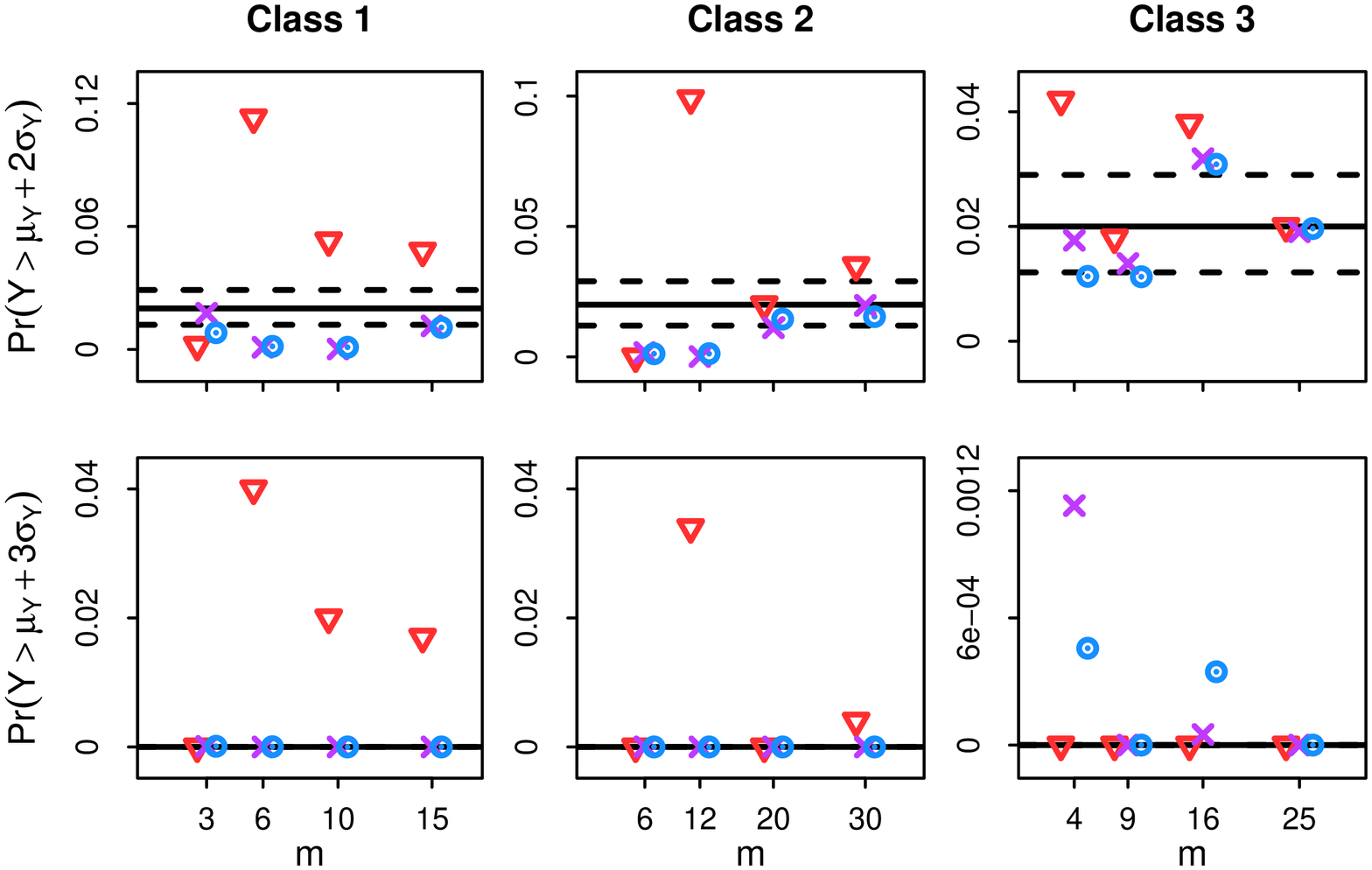}
\end{figure}

\begin{figure}
\centering
\caption{Probability density function (PDF) estimation validation results for the two-dimensional function in \eqref{nonlinear}. Design classes are in columns and design size used to build the surrogates increases further down the rows. Polynomial chaos, squared exponential and Mat\'{e}rn Gaussian process surrogate model PDF estimates are shown as red, purple and blue lines respectively, while the simulator output PDF is shown as a black line. Gaussian process emulators also have $95\%$ confidence intervals about their estimates.}
\label{exampledensity}
\epsfysize=110mm \epsfbox{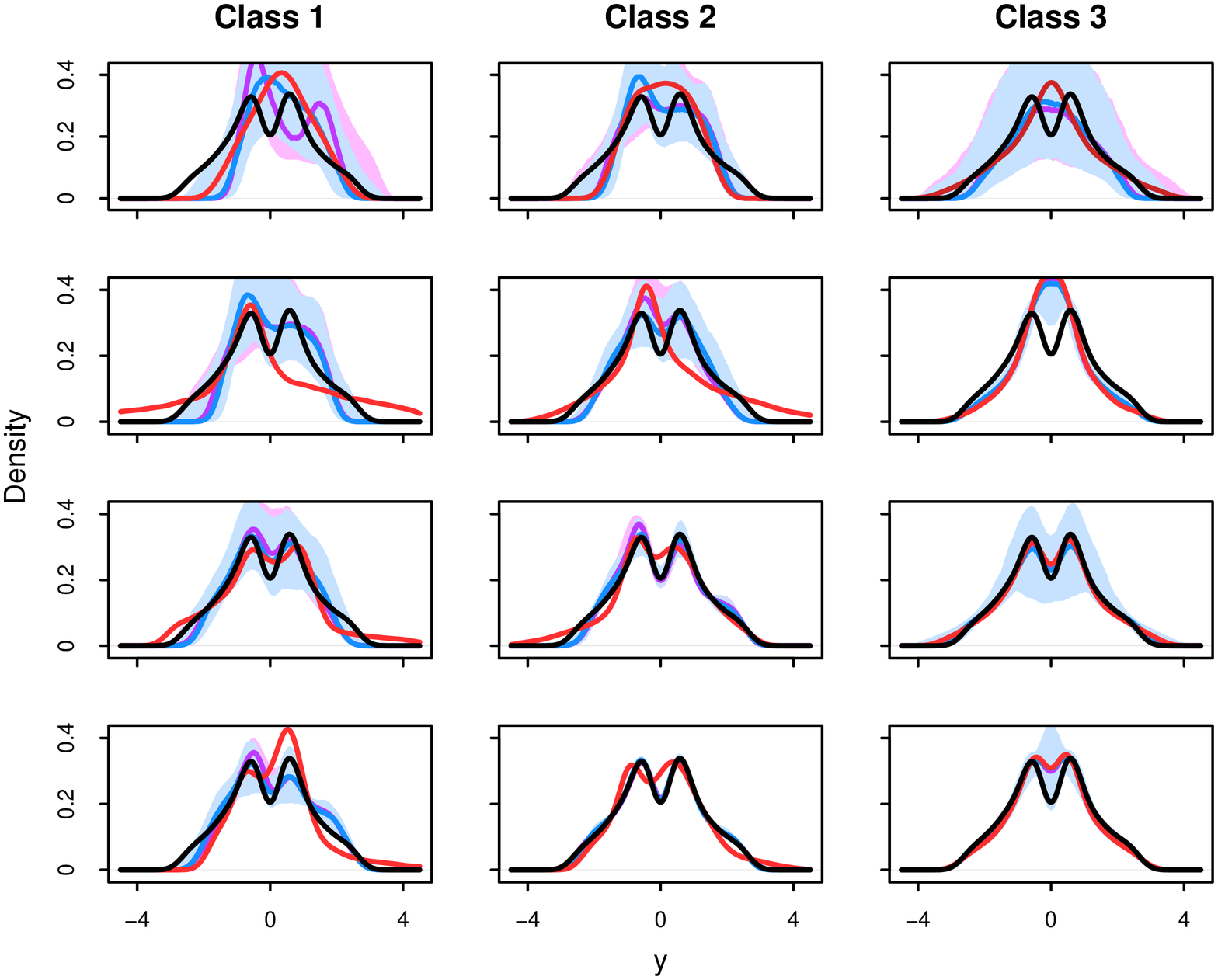}
\end{figure}

This toy simulator, plotted in Figure \ref{exampleplot}, is clearly a nonlinear function of the input parameters. As a result of this, we would expect a low-order polynomial expansion to perform poorly. We repeat our experiments for the toy simulator by fitting PC, SE GP and M GP surrogates to class 1, 2 and 3 designs and evaluating validation metrics based on an independent 1000-point Latin hypercube design. Once again, the design sizes are restricted by the polynomial chaos methods chosen, as well as the input dimension ($n=2$). Validation results for the two-dimensional function are shown in Figures \ref{example}, \ref{exampleexceed} and \ref{exampledensity}. Interestingly, polynomial chaos accuracy for this function is not so poor for class 1 designs as it was for the \verb+adJULES+ and \verb+VEGACONTROL+ experiments. However, both Gaussian process emulators still give more accurate and stable results so they are preferred here. One of the most striking results is the difference between the two methods for class 2 designs. For all validation metrics and the majority of design sizes, Gaussian process emulators consistently perform better than polynomial chaos surrogates. We also note for class 3 designs (tensor grids where polynomial chaos surrogates where slightly preferred), the performance gap between the two methods has closed considerably. Finally, this nonlinear function exhibits a bimodal PDF that all surrogate methods find difficult to estimate, especially for small design sizes. Nonetheless, a preference for GP methods is still apparent. These results provide good evidence that Gaussian process emulators are more suitable for modelling nonlinear simulator behaviour. Nevertheless, polynomial chaos surrogates still provide respectable accuracy for a fraction of the cost. \\

Finally we reflect on the practicality of polynomial chaos and Gaussian process surrogates --- once built, what do they give us? It is clear that both methods offer a fast-approximation to the simulator at any untried input configuration, which is useful in performing Monte Carlo based uncertainty quantification tasks such as uncertainty or sensitivity analysis. Crucially though, Gaussian process emulators not only provide a single prediction but have the advantage of readily available uncertainty information due to the distributional assumptions. This not only allows the user to see where the emulator is most uncertain in the design space, but uncertainty can be fed through to various validation metrics as we showed in our examples. This is not the case for polynomial chaos, and thus some preference may be found for Gaussian process emulators if uncertainty information is required by the practitioner. \\

In conclusion, polynomial chaos and Gaussian process emulation are both state-of-the-art approaches with the focus of performing surrogate-based uncertainty quantification for computationally expensive simulators. Despite this common goal, their respective communities tend to work on their own problems and are only just starting to recognise each other. We hope the work presented here goes some way towards bringing the two communities closer together, or at least is useful advice for practitioners in the field who may be unsure of what method to use. Scope for future work in the area is vast, including similar comparisons for higher dimensional problems and more rigorous or theoretical assessments of the two methods in terms of flexibility, practicality and computational cost. In this work we were primarily concerned with the prediction objective, that is, using a surrogate to give an estimate of the simulator output at an untried input location. Comparisons of the two methods for other specific uncertainty quantification tasks, such as calibration or sensitivity analysis, would also be valuable. Development of a hybrid method combining polynomial chaos and Gaussian process emulation may also be possible, and such an approach could draw upon the advantages of both approaches to give a more accurate surrogate.  

\clearpage
\bibliography{Nathan}
\bibliographystyle{ametsoc}

\end{document}